\documentstyle[12pt,amssymb]{amsart}

\newcommand\res{ \restriction }

\newcommand\OR{{\mbox{OR}}}

\newcommand\ult{{\mbox{ult}}}
\newcommand\pred{{\mbox{pred} \,}}
\newcommand\rt{{\mbox{root} \,}}
\newcommand\decap{{\mbox{decap}}}

\newcommand\kvec{{\vec k}}

\newcommand\Evec{{\vec E}}

\newcommand\gbvec{{{\vec {\gb} }}}
\newcommand\gkvec{{{\vec {\gk} }}}
\newcommand\etavec{{{\vec {\eta} }}}
\newcommand\glvec{{{\vec {\gl} }}}

\newcommand\gsvec{{{\vec {\gs} }}}

\newcommand\gLvec{{{\vec {\gL} }}}

\newcommand\cNvec{{{\vec {\cal N} }}}
\newcommand\cPvec{{{\vec {\cal P} }}}
\newcommand\cQvec{{{\vec {\cal Q} }}}
\newcommand\cRvec{{{\vec {\cal R} }}}
\newcommand\cSvec{{{\vec {\cal S} }}}

\newcommand\cQtil{{\widetilde {\cQ}}}
\newcommand\cTtil{{\widetilde {\cT}}}

\newcommand\gbtil{{\widetilde {\gb}}}

\newcommand\gptil{{\widetilde {\gp}}}
\newcommand\gdtil{{\widetilde {\gd}}}

\newcommand\etatil{{\widetilde {\eta}}}

\newcommand\sdot{{\dot s}}
\newcommand\Edot{{\dot E}}
\newcommand\Fdot{{\dot F}}
\newcommand\ggdot{{\dot \gg}}
\newcommand\gmdot{{\dot \gm}}
\newcommand\gndot{{\dot \gn}}
\newcommand\gtdot{{\dot \gt}}

\newcommand\Ebar{{\overline E}}

\newcommand\Tbar{{\overline T}}
\newcommand\Wbar{{\overline W}}

\newcommand\cTbar{{\overline {\cT}}}

\newcommand\gkbar{{\overline \gk}}
\newcommand\glbar{{\overline \gl}}

\newcommand\gObar{{\overline \gO}}
\newcommand\gUbar{{\overline \gU}}

\newcommand\cD{{\cal D}}

\newcommand\cG{{\cal G}}
\newcommand\cH{{\cal H}}

\newcommand\cJ{{\cal J}}

\newcommand\cL{{\cal L}}
\newcommand\cM{{\cal M}}
\newcommand\cN{{\cal N}}

\newcommand\cP{{\cal P}}
\newcommand\cQ{{\cal Q}}
\newcommand\cR{{\cal R}}
\newcommand\cS{{\cal S}}
\newcommand\cT{{\cal T}}
\newcommand\cU{{\cal U}}
\newcommand\cV{{\cal V}}

\newcommand\ra{\rightarrow}

\newcommand\lra{\longrightarrow}

\newcommand\Lra{\Longrightarrow}

\newcommand\Llra{\Longleftrightarrow}

\newcommand{\ga}{\alpha}     
\newcommand{\gb}{\beta}      
\renewcommand{\gg}{\gamma}   
\newcommand{\gd}{\delta}     
\renewcommand\ge{\varepsilon}
\newcommand{\gz}{\zeta}

\newcommand{\gk}{\kappa}  
\newcommand{\gl}{\lambda}    
\newcommand{\gm}{\mu}        
\newcommand{\gn}{\nu}

\newcommand{\gp}{\pi}        
\newcommand{\gr}{\rho}       
\newcommand{\gs}{\sigma}     
\newcommand{\gt}{\tau}

\newcommand{\go}{\omega}

\newcommand{\gG}{\Gamma}

\newcommand{\gL}{\Lambda}

\newcommand{\gS}{\Sigma}     
       
\newcommand{\gU}{\Upsilon}

\newcommand{\gO}{\Omega}

\newcommand{\ha}{\aleph}

\renewcommand{\models}{\vDash}

\newcommand{\dom}{\mbox{dom}}
\newcommand{\ran}{\mbox{ran}}
\newcommand{\crit}{\mbox{crit}}
\newcommand{\card}{\mbox{card}}

\newcommand{\cf}{\mbox{cf}}
\newcommand{\lh}{\mbox{lh}}

\newcommand{\id}{\mbox{id}}

\begin{document}
\baselineskip 4ex 

\title{The covering lemma up to a Woodin cardinal}{}

\author{W.J. Mitchell, E. Schimmerling, and J.R. Steel}
\thanks{The research of the first author was partially supported by
National Science Foundation grants DMS-9001027 and DMS-9306286.}
\thanks{The research of the second author was partially supported by a 
National Science Foundation Mathematical Sciences Postdoctoral Research
Fellowship.}
\thanks{The research of the third author was partially supported by
National Science Foundation grant DMS 92-06946.}

\maketitle

\vskip 6ex
\noindent
{\sc {\S}1 \ \underline{Introduction}}
\vskip 4ex

We say that a cardinal $\gd$ is countably closed if
$\gg^{\ha_0} < \gd$ whenever $\gg < \gd$.
In this paper, we prove the following theorem.

\vskip 3ex
\noindent
{\sc Theorem 1.1.}
Suppose that $\gO$ is a measurable cardinal,
and that there is no inner model with a Woodin cardinal.
Let $K$ be the core model constructed in $V_\gO$.
Suppose that $\gk < \gO$ is a $K$-cardinal
such that $\card (\gk )$ is countably closed.
Let $\gl = ( \gk^+ )^K$.
Then $$\gl < \gk^+ \ \Lra \ \cf (\gl ) = \card ( \gk ).$$
In particular, $K$ computes successors of countably closed, singular
cardinals correctly.
\vskip 3ex

In the theorem, $K$ is the core model for one Woodin cardinal, 
as introduced in Steel \cite{St1}.

In order to trace the history behind the theorem,
it will be convenient to
make the following definitions.
Consider an inner model $M$ and an ordinal $\ga$.
Let us say that CP($M$,$\ga$) holds whenever
$$\forall X \subseteq \ga \ \exists Y \in M \ 
(X \subseteq Y \ \& \ \card (Y) \leq \card (X) + \ha_1 )$$
and that WCP($M$,$\ga$) holds whenever
$( \ga^+ )^M = \ga^+$.
The theorem concludes that if $\gk$ is a singular, countably closed
cardinal, then WCP($K$,$\gk$) holds.

Jensen's Covering Lemma for $L$ says that either
CP($L$,$\ga$) holds for every ordinal $\ga$,
or else $0^\#$ exists; see \cite{De}.
Dodd and Jensen introduced their core model in \cite{DoJe1} and 
found an analogue of the result for $L$ in \cite{DoJe2}.
A consequence of their work is that either
CP($K$,$\ga$) holds for every ordinal $\ga$,
or there is a measurable cardinal in $K$,
or $0^\dagger$ exists and is an element of $K$.
The Dodd-Jensen result is extremely useful in proving that
various set theoretic hypotheses have large cardinal consistency strength of
at least one measurable cardinal:
one merely has to violate CP($K$,$\ga$).
However, this method for obtaining lower bounds on consistency strength
does not extend, at least without modification, to larger core models.
This is because if $\gm$ is measurable,
then CP($K$,$\gm$) fails in any Prikry generic forcing extension.
Luckily,
many of the applications of the Covering Property 
are really consequences of WCP($K$,$\gk$) for some singular cardinal $\gk$.
It is easy to see that if $\gk$ is a singular cardinal,
then CP($M$,$(\gk^+ )^M$) implies WCP($M$,$\gk$).  The weak covering property
that $K$ computes successors of singular cardinals correctly, is expected
to be true for much larger core models.

Indeed, above one measurable cardinal, weak covering properties are
more a part of the core model theory, than a consequence thereof.
In Mitchell \cite{Mi1,Mi2}, where the working assumption is that there is no
inner model with a measurable cardinal $\gm$ of Mitchell order $\gm^{++}$,
Mitchell defines the countably complete core model, $K^{cc}$, 
and shows that for all
countably closed, singular cardinals $\gk$, WCP($K^{cc}$,$\gk$) holds.
Under the same working assumption, Mitchell goes on 
to define $K$ as a transitive collapse of a certain elementary 
substructure of $K^{cc}$ and to develop the basic properties of $K$
using the iterability and weak covering properties of $K^{cc}$.

The working assumption in Steel \cite{St1} is only that there is no inner model
with a Woodin cardinal.  Steel defines a preliminary model, $K^c$,
the background certified core model, and shows that $K^c$ is iterable in
a sense involving iteration trees.  But in order to establish
a form of the weak covering property for $K^c$, Steel makes the {\it ad hoc}
assumption that $K^c$ is constructed inside some $V_\gO$ with $\gO$
a measurable cardinal in $V$.  He shows that
WCP($K^c$,$\ga$) holds for almost every $\ga < \gO$,
in the sense of any normal measure on $\gO$,
and goes on to use this ``cheapo'' covering property
and the iterability of $K^c$ to isolate
and develop the basic properties of $K$.  Our Theorem 1.1 implies that
$K$ satisfies the Mitchell weak covering property in this setting.

Some important open questions regarding Theorem 1.1 remain.
Can the assumption that there is a measurable cardinal
be dropped?
\footnote{See Schimmerling \cite{Sch3} for some partial results.}
What about the assumption that
$\card ( \gk )$ is countably closed, 
can it be replaced by just $\gk \geq \ha_2$?
\footnote{Mitchell and Schimmerling \cite{MiSch} have since answered this
question positively,
by a proof that builds directly on the results presented here.}
Towards a positive answer to these questions,
Jensen showed in \cite{Je2} that the assumption that there is a
measurable cardinal above $\gk$, and that $\gk$ is countably closed, 
are not needed if, instead,
one assumes that $0^{\P}$ does not exist and $\gk \geq \ha_3$.
Another question is whether Theorem 1.1 can be extended to richer core models.
The proof of Theorem 1.1, 
combined with the theories developed in
Steel \cite{St3, St4}, Schimmerling-Steel \cite{SchSt}, and Woodin
\cite{W},
yields extensions of Theorem 1.1 to core models for mice with many Woodin
cardinals.

This paper is organized as follows.
Building on Mitchell-Steel \cite{MiSt}, Steel \cite{St1}, and 
Schimmerling \cite{Sch},
we develop further the fine structure of $K$ in {\S}2.
While we shall make some effort to remind the
reader of the relevant definitions and facts,
we shall, in the end, assume
that the reader is quite familiar with these earlier works.
The proof of Theorem 1.1 is given in {\S}3.
We encourage the reader to look ahead to {\S}3 for additional
motivation for the technical results in {\S}2.

In proofs of covering properties for smaller
core models, it is shown that for particular embeddings
$\pi : N \lra V$, if we compare $K^N$ with $K$
and let $\cR$ be the ultrapower of the final model on
the $K$-side by the long extender from $\pi$,
then $\cR$ is iterable in a way that guarantees
that $\cR$ is a level of $K$.
The main new idea in the proof of Theorem 1.1
is to prove the required iterability, and more, by induction
on the models occurring on the $K$-side of the comparison.
The strengthening of the iterability hypothesis keeps the
induction going.

The idea of proving more than the required iterability by induction,
is due to Mitchell,
and was first described by him in a note dated October, 1990.
All three authors contributed to the work of converting that 
note into this paper.   
In the end, however, we were not able to show that the structure $\cR$
of the previous paragraph is a level of $K$.  As far as we know,
$\cR$ may itself be class size, or worse, $\cR$ may fail to be 
a potential premouse, making it impossible that $\cR$
be a level of $K$.
The possibility that $\cR$ is a proper class
requires rearranging the proof in a way that is not
particularly more complicated and leads to an interesting
new perspective.
Our method for dealing with the possibility that $\cR$ is not a 
premouse, is due to Schimmerling, and was inspired
by \cite{Sch}.

\vskip 6ex
\noindent
{\sc {\S}2 \ \underline{A Barrage of Fine Structure}}
\vskip 4ex

In this section we extend the fine structure established
in \cite{MiSt}, \cite{St1}, and \cite{Sch}, and prove some
facts that will be used in the proof of Theorem 1.1.
We expect that the reader is quite familiar with those works,
but we will try, whenever possible, to remind the reader
of what things mean.

\vskip 6ex
\noindent
{\sc {\S}2.1 \ \underline{The Dodd-squash}}
\vskip 4ex

Suppose that $F$ is a $(\gk , \gb)$-extender (or pre-extender)
over a potential premouse (ppm) $\cM$.  $\gk$ is the {\bf critical
point of} $F$ and $\gb$ is the {\bf length of} $F$;
we write $\crit (F) = \gk$ and $\lh (F) = \gb$.

Suppose that $\xi \geq \gk$ and $t \in [ \gb ]^{< \go}$.
We say that $\xi$ is a {\bf generator of} $F$
{\bf relative to} $t$ iff 
$\xi \not= [ a \cup t , f ]^\cM_F$ for any $a\in [ \xi ]^{<\go}$
and $f \in |\cM |$ with $f : [\gk ]^{|a \cup t|} \lra \gk$.
We say that $\xi$ is {\bf generator of} $F$ iff $\xi$
is a generator of $F$ relative to $\emptyset$.

Note that whether or not $\xi$ is a generator of $F$ relative to $t$ depends
on $\xi$, $t$, $F$, and $\cP (\gk ) \cap | \cM |$.
But since $F$ determines $\cP (\gk ) \cap | \cM |$,
the dependence on $\cM$ is
merely nominal.  This will be true of many of the properties 
of $F$ which we shall discuss below.

By the {\bf strict supremum of the generators of} $F$ we mean
$$\gn (F) =
\sup ( \{ \ \xi +1 \ | \ \xi \mbox{ is a generator of } F \ \} 
		\cup (\gk^+ )^\cM ).$$
Let $\cP = \ult_0 (\cM , F)$
and let $j : \cM \lra \cP$ be the ultrapower map.
By the {\bf trivial completion of} $F$, we mean the 
$(\gk , ( \gn ( F )^+ )^\cP)$-extender derived from $j$.
Recall that if $\Evec$ is an extender from a good extender sequence,
then every $E_\ga$ is its own trivial completion.
Assume for the rest of this subsection 
that the trivial completion of $F$ is $F$ again;
in particular, $\gb = ( \gn ( F )^+ )^\cP$.
By $F^*$ we mean the 4-ary 
relation given by:
\begin{tabbing}
\ \ \ \ \ \ \ $(\xi , \gg , a , x) \in F^*$ if and only if
\= $\gk < \xi < ( \gk^+ )^\cM$,\\
\> $\gg < \gb$,\\
\> $F \ \cap 
	\ ( [\gn (F) ]^{ < \go } \times J^\cM_\xi ) \in J^\cP_\gg$, and\\
\> $(a, x) \in 
	F \ \cap \ ( [\gg ]^{ < \go } \times J^\cM_\xi ).$
\end{tabbing}
We purposely allow the possibility that $(\gk^+ )^\cM < (\gk^+ )^\cP$.
$F^*$ codes $F$ by the fragments of $F$.

Varying slightly from \cite{MiSt},
our convention on ppm's is that if $\cM$ is an active ppm,
then
$\Fdot^\cM = F^*$ for some extender $F$ that is its own trivial
completion.  We call $F$ (or $F^*$) 
the {\bf last extender of $\cM$} in this case.
The advantage is that ppm's defined this way are amenable structures,
and there is no need to consider {\it restricted $\gS_n$}
formulas.
The equivalence between this approach and the official approach
taken in \cite{MiSt} (which we avoid going into here)
is explained in \cite{MiSt}.

The following definitions have their roots in some unpublished
notes of A.J. Dodd.  The reader should consult Section 3 of \cite{Sch}.

We define the {\bf Dodd-parameter of} $F$ = $s(F)$ by induction.
$s(F)(0)$ is the last generator of $F$ above $(\gk^+ )^\cM$,
should it exist, and is left undefined otherwise.
If $s(F)(0)$, ... , $s(F)(i)$ are defined,
then $s(F)(i+1)$ is the last generator of $F$ relative to 
$\{ s(F)(0)$, ... , $s(F)(i) \}$ that is above $(\gk^+ )^\cM$,
should it exist, and is left undefined otherwise.
Since $s(F)(i+1) < s(F)(i)$, we end up with 
$s(F) \in 
[ \ga - (\gk^+ )^\cM ]^{< \go}$.

By the {\bf Dodd-projectum} of $F$ we mean the ordinal
$\gt (F)$ = the strict supremum of the generators of $F$
relative to $s(F)$, that is:
$$ \gt (F) = \sup ( \{ \ \xi +1 \ | \ \xi \mbox{ is a generator of } F 
		\mbox{ relative to } s(F) \ \}
		\cup (\gk^+ )^\cM ).$$
The {\bf Dodd-squash} of $F$ is
$F \res ( \gt (F) \cup s(F) )$,
usually viewed as coded by a subset of
$$[ \gt \cup \{ {\check s_0}, {\check s_1}, \dots \} ]^{< \go}
\times J^\cM_{(\gk^+)^\cM} \ ,$$
where each $\check s_i$ is a constant symbol in HF (the collection
of hereditarily finite sets) representing $s_i$.

Suppose that $\cM$ is an active ppm with $F^* = \Fdot^\cM$ = the 
last extender of $\cM$.
By the {\bf Dodd-squash} of $\cM$ we mean the structure
$$ \cM^{Ds} = 
\langle \ J^\cM_{\gt (F)} \ , \ \in \ , \  \Edot^\cM \res \gt (F) \ , 
\ F \res (\gt (F) \cup s(F)) \ \rangle.$$
In this situation, we also write
$\gtdot^\cM = \gt (F)$ and $\sdot^\cM = s(F)$.
In \cite{MiSt}, active ppm are divided up into three types:
\begin{list}{}{}
\item[ ]{
$\cM$ is type I iff $F$ has only $\gmdot^\cM = \crit (F )$ as a generator}
\item[ ]{
$\cM$ is type II iff $F$ has a last generator, but is not type I}
\item[ ]{
$\cM$ is type III iff the generators of $F$ have limit order type}
\end{list}
Just to remind the reader:
\begin{list}{}{}
\item[ ]{If $\cM$ is type I, then
$\gtdot^\cM = ( \gmdot^+ )^\cM = ( \crit (F)^+ )^\cM = \gndot^\cM$
and $\sdot^\cM = \emptyset$.}
\item[ ]{
If $\cM$ is type II,
then $\gtdot^\cM < \gndot^\cM = \gn (F) = s(F)(0)+1$.}
\item[ ]{If $\cM$ is type III, then 
$\gtdot^\cM = \gndot^\cM = \gn (F)$, $\sdot^\cM = \emptyset$,
and $\cM^{Ds} = \cM^{sq}$.}
\end{list}

Some basic facts about the Dodd-squash are established in Section 3 of 
\cite{Sch}.  Most notably, a proof of the following theorem of
Steel is given.

\vskip 3ex
\noindent
{\sc Theorem 2.1.1.}
Suppose that $\cM$ is a 1-sound, 1-small, 0-iterable premouse.
Let $F = \Fdot^\cM$, $s = \{ s_0 , \dots , s_k \} = \sdot^\cM$,
and $\gt = \gtdot^\cM$.  Then:
\begin{list}{}{}
\item[(A)]{ $\gt$ is a cardinal in $\cM$.}
\item[(B)]{For $i\leq k$,
$F_i = F \res (s_i \cup (s \res i )) \in | \cM |$.}
\item[(C)]{$\cM^{Ds}$ is amenable.}
\end{list}
\vskip 3ex

The set $F_i$ in (B) we call the $i^{th}$ {\bf Dodd-witness for} $F$
(or $s(F)$ or $\cM$ or $\cM^{Ds}$),
and (B) asserts that $F$ (or 
$s(F)$ or $\cM$ or $\cM^{Ds}$) is {\bf Dodd-solid}.

The next lemma extends Lemma 4.4 of \cite{Sch}.

\vskip 3ex
\noindent
{\sc Lemma 2.1.2.}
Suppose that $\cM$ is a type II ppm.  Let $\gt = \gtdot^\cM$,
$s =  \{ s_0 ,\dots , s_k \} = \sdot^\cM$, $F = \Fdot^\cM$, 
$\mu = \gmdot^\cM = \crit (F)$.
Suppose that $\cM^{Ds}$ is amenable.
Also, suppose that for some $\gk$ such that
$$(\gm^+ )^\cM \leq \gk < \gt$$
and some parameter $q$ from $\cM$,
$$|\cM | = H^\cM_1 (\gk \cup q ).$$
Then there is a parameter $x \in [ \gt ]^{< \go}$ such that
$$\sup ( \gt \cap H^{\cM^{Ds}}_1 (x) ) = \gt.$$
\vskip 3ex

\noindent
{\sc Proof.}
Choose $x \in  [ \gt ]^{< \go}$ with the following properties:
\begin{list}{}{}
\item[1.]{$\gk \in x$}
\item[2.]{$s \in H^\cM_1 (x \cup q )$}
\item[3.]{There is some $f \in J^\cM_{(\gm^+)^\cM}$
	such that 
	$[ x \cup s , f ]^\cM_F = \langle q , \ggdot^\cM \rangle$\\
(recall that $\ggdot^\cM$ codes the last initial segment of
$F$).}
\end{list}
Let $\gs = \sup ( \gt \cap H^{\cM^{Ds}}_1 (x) )$.
By 1, $\gs > \gk \geq (\mu^+ )^\cM$.
Put $$\cQtil = \cM^{Ds} \res \gs
= \langle \ J^\cM_\gs \ , \ \in \ , \ \Edot^\cM \res \gs \ , \ F \res (\gs \cup s )
\ \rangle.$$
Let $\cQ$ be the {\bf de-Dodd-squash of} $\cQtil$,
that is, the structure
$$\langle \ J^\cP_\gb \ , \ \in \ , \  \Edot^\cP \res \gb \ , \ G^* \ \rangle,$$
where $\cP = \ult_0 (\cQtil , F \res (\gs \cup s ))$,
$$\gb = 
\left[ \{ s_0 \} , 
\{ \xi \} \mapsto (\xi^+)^\cQtil \right]^\cQtil_{F\res (\gs \cup s )} \ ,$$
and $G$ is the $(\mu , \gb )$-extender derived from the ultrapower
embedding $\cQtil \lra \cP$.

Let $\pi : \cQ \lra \cM$ be the canonical embedding.
A straightforward calculation shows that 
$\pi$ is a 0-embedding of $\cQ$ into $\cM$
with $\crit (\pi ) \geq \gs $.
Note that 
$$| \cQ | = H^\cQ_1 (\gs \ \cup \ \pi^{-1} (s)).$$
Using 2 and 3, we see that 
$$| \cQ | = H^\cQ_1 (\gs \ \cup \ \pi^{-1} (q)).$$
But since 
$$| \cM | = H^\cM_1 (\gs \ \cup \ q),$$
$\pi = \id$
and $\gs = \gt$.
\qed {\tiny \ \ Lemma 2.1.2}
\vskip 3ex

From the lemma, standard fine structural calculations show:

\vskip 3ex
\noindent
{\sc Corollary 2.1.3.}
Under the same hypothesis,
the ordinals $\gt$, $(\gk^+ )^\cM$, $\OR^\cM$, and $(\mu^+ )^\cM$
have countable cofinality, definably over $\cM$.
\vskip 3ex

The next lemmas tell us how $\gtdot$ and $\sdot$ change as we iterate.

\vskip 3ex
\noindent
{\sc Lemma 2.1.4.}
Suppose that $E$ is a $(\gk , \gb )$-extender over a ppm $\cM$,
and $\gk < \gtdot^\cM$.
Suppose that $\cM^{Ds}$ is amenable and Dodd-solid;
say $\gt = \gtdot^\cM$ and $s = \{ s_0 \, , \dots , \, s_l \} = \sdot^\cM$.
Let $i:\cM \lra \cP = \ult_0 ( \cM , E)$ be the ultrapower
map.
Then $\sdot^\cP = i(s)$,
$\gtdot^\cP = \sup ( i '' \gt )$,
and $\cP^{Ds}$ is amenable.
\vskip 3ex

\noindent
{\sc Proof.} The proof is essentially the same as that of
Lemma 9.1 of \cite{MiSt}.
We have the commutative diagram:
$$\begin{array}{ccc}
\ult_0 (\cM , \Fdot^\cM ) & 
\stackrel{j}{\lra} & \ult_0 (\cP , \Fdot^\cP ) \cr
\gp \uparrow & \ & \uparrow \gs \cr
\cM & \stackrel{i}{\lra} & \cP = \ult_0 (\cM , E )
\end{array}$$
with $i = j \res |\cM |$.
Since $i$ is a 0-embedding and $\cM$ is Dodd-solid,
$$i( \Fdot^\cM \res (s_n \cup (s|n)) ) = 
\Fdot^\cP \res (i (s_n ) \cup i(s|n) ) \in | \cP |.$$
Since $i$ is a 0-embedding and $\cM^{Ds}$ is amenable,
for each $\xi < \gt$, $$i ( \Fdot^\cM \res (\xi \cup s ) )
= \Fdot^\cP \res (i ( \xi ) \cup i(s) ) \in 
|\cP |,$$ and
$\gt^* = \sup ( i '' \gt )$ is a limit of generators of $\Fdot^\cP$
relative to $i(s)$.
So it is enough to see that no $\eta \geq \gt^*$ is a generator of
$\Fdot^\cP$ relative to $i(s)$.

So suppose that $\gt^* \leq \eta < \OR^\cP$.
Then $\eta = j(f)(b)$ for some $f \in |\cM |$
and $b \in [ \gb ]^{< \go}$ with $f: [ \gk ]^{|b|} \lra \gk$.
Since $f \in |\cM | \subset |\ult_0 (\cM , \Fdot^\cM ) |$,
$f = \pi ( g )(a \cup s)$ for some $g\in |\cM |$ and 
$a \in [ \gt ]^{< \go }$ with 
$g : [ \gmdot^\cM ]^{|a \cup s|} \lra |\cM |$.
Thus
$$\eta = j(f)(b)=j(\pi (g )(a \cup s ))(b)
= j(\pi (g ))(j(a) \cup j(s) )(b)
= \gs ( i (g))(i(a) \cup i(s))(b).$$
So
$\eta = [i(a) \cup b \cup i(s), h ]^\cP_{\Fdot^\cP}$
where
$$h(u) = i(g)(u^{i(a) \cup i(s) \ , \ i(a) \cup b \cup i(s)})(u^{b\ ,\ i(a) \cup b
\cup i(s)}).$$
Since $i(a) \in [\gt^* ]^{< \go}$ and $b \in [ \gb ]^{< \go } \subseteq
[ \gt^* ]^{< \go } $, we are done proving Lemma 2.1.4.\\
\qed {\tiny \ \ Lemma 2.1.4}
\vskip 3ex

We say that $\cM$ is {\bf Dodd-solid above $\gk$}
if the $i$'th Dodd-solidity witness for $\Fdot^\cM$
is an element of $|\cM |$ whenever $\sdot^\cM_i \geq \gk$.
An argument similar to that just given shows:

\vskip 3ex
\noindent
{\sc Lemma 2.1.5.}
Suppose $E$ is a  $(\gk , \gb)$-extender over $\cM$ with $\gtdot^\cM \leq
\gk$.
Suppose that
$\cM$ is Dodd-solid above $\gk$.
Let $j:\cM \lra \cP = \ult_0 (\cM , G)$ be the ultrapower map.
Then $\gtdot^\cP \leq \gn (E)$,
$$j\left( \{ \sdot^\cM_i | \sdot^\cM_i \geq \gk \} \right)
= \{ \sdot^\cP_i \mid \sdot^\cP_i \geq \gn (E) \} \ ,$$
and $\cP$ is Dodd-solid above $\gn (E)$.
\vskip 3ex

Note that under the hypothesis of Lemma 2.1.5, it is possible that $\cP^{Ds}$
is not amenable or that $\sdot^\cP$ is not Dodd-solid; see 
Section 3 of \cite{Sch} for an examples using $0^{\P}$.
The main corollary to the previous lemmas that we shall need is:

\vskip 3ex
\noindent
{\sc Corollary 2.1.6.}
Suppose that $\cT$ is a padded $\go$-maximal iteration tree on a
1-small, $\gO$-iterable weasel $W$ with models $W_\eta$ for $\eta < \lh
(\cT )$.
Suppose that $\eta < \lh ( \cT )$ and $\gk$ is a cardinal in $W_\eta$
such that 
$$\gk \geq \sup \left( \{ \nu^\cT_\gb \ | \ \gb + 1 \leq_T \eta \} \right)$$
(recall that $\nu^\cT_\gb = \gn (E^\cT_\gb )$ for such iteration trees).
Suppose that $\cP$ is an initial 
segment of $W_\eta$ such that $\gr_1^\cP \leq \gk$ and
$\gk \geq (\gmdot^+)^\cP$.
Then either $\gtdot^\cP \leq \gk$ and $\cP$ is Dodd-solid above $\gk$,
or else the ordinals $\gtdot^\cP$, $(\gk^+ )^\cP$,
$\OR^\cP$, and $(\gmdot^+ )^\cP$ all have countable cofinality.

\vskip 6ex
\noindent
{\sc {\S}2.2 \ \underline{For Weasels Only}}
\vskip 4ex

By a {\bf weasel}, we mean a premouse of height $\gO$.
As we have seen in \cite{Mi1}, \cite{Mi2}, and \cite{St1},
the hull and definability properties for weasels play a role
similar to the role that soundness plays for set premice.
We take this analogy a bit further in this subsection.

\vskip 3ex
\noindent
{\sc Definition 2.2.1.}
Let $A_0$ be the class defined in Section 1 of \cite{St1}.
Suppose that $W$ be a weasel,
$s \in [\gO ]^{< \go}$, and  $\gk < \gO$.
We say that $W$ has the {\bf $s$-hull property at $\gk$} if
and only if 
$$\cP ( \gk ) \subseteq \mbox{ the transitive 
collapse of } H^W ( \gk \cup s \cup \gG )$$
whenever $\gG$ is $A_0$-thick in $W$.
We say that $W$ has the {\bf $s$-definability property at $\gk$} if
and only if 
$$\gk \in H^W (\gk \cup s \cup \gG )$$
whenever $\gG$ is $A_0$-thick in $W$.
\vskip 3ex

In the definitions that we are about to give, there is a somewhat
suppressed
dependence on our choice of $\gU < \gO$.
In the proof of Theorem 1.1, $\gU$ will be a cardinal $> \gl$.

By analogy with both the standard parameter and the Dodd-parameter,
we define the {\bf class-parameter} $c(W )$ by induction.
Let $c(W)(0)$ be the largest $\xi < \gU   $ such that $W$
does not have the definability property at $\xi$; if no such $\xi$
exists, then leave $c(W)(0)$ undefined.
Suppose that we have defined $c(W )(0) , \dots ,  c(W )(i)$.
Let $c(W  )(i+1)$ be the largest $\xi < \gU$ 
such that
$W$ does not have the 
$\{ c(W )(0) , \dots , c(W )(i) \}$-definability property at $\xi$;
if no such $\xi$ exists, then leave $c(W  )(i+1)$ undefined.
Of course, $c(W )(i+1) < c(W )(i)$, so we are in fact getting 
$c(W ) \in [\gU    ]^{<\go}$.
We define the {\bf class projectum} of $W$ to be 
$\gk ( W  )$ = the supremum of the ordinals $\xi$ such that 
$W$ does not have the $c(W )$-definability property at $\xi$.
Notice that if $\gG$ is $A_0$-thick in $W$, then
we have the following analog to soundness:
$$H^W ( \gk (W  ) 
\, \cup \, c(W ) \, \cup \, \gG ) \ \cap \ J^W_{\gU   }
= J^W_{\gU   }.$$
In particular, $W$ has the $c(W )$-hull property at $\gk (W )$
when $\gk (W) < \gU$.

\vskip 8ex
\noindent
{\sc {\S}2.3 \ \underline{Protomice}}
\vskip 4ex

\noindent
{\sc Definition 2.3.1.}
Suppose that $\cR$ and $\cM$ are a ppm's,
$F$ is a $(\gk , \gb )$-extender over $\cM$,
$(\gk^+ )^\cM \leq ( \gk^+ )^\cR$,
and $\cM$ and $\cR$ agree below $(\gk^+ )^\cM$.
Then we say that 
$F$ is a {\bf $(\gk , \cM , \gb )$-extender fragment over} $\cR$.
\vskip 3ex

So extenders are extender fragments, but not necessarily the other
way around.
Note that if $F$ is a $(\gk , \cM , \gb )$-extender fragment over $\cR$,
then $F$ is a $(\gk , \cN , \gb )$-extender fragment over $\cR$
for any initial segment $\cN$ of $\cR$ with 
$(\gk^+ )^\cN 
= (\gk^+ )^\cM$. 
Extender fragments arise naturally in the proof of weak
square given in \cite{Sch}.
They are a fine structural analogue of the background certificates
used to build $K^c$ in \cite{St1}.
In the proof of Theorem 1.1,
we shall often encounter structures that look exactly like a premice,
except that their last predicates code extender fragments
that are not total extenders.

\vskip 3ex
\noindent
{\sc Definition 2.3.2.}
A {\bf protomouse} is a structure $\cR$ in the language of ppm's
which is of the form
$$\langle \ J^\Evec_\ga \ , \ \in \ , \ \Evec \res \ga \ , \ F^* \
\rangle,$$
where
$$\cR' = 
\langle \ J^\Evec_\ga \ , \ \in \ , \ \Evec \res \ga \ \rangle$$
is a passive premouse,
and $F$ is a 
$(\gm , \cM , \ga )$-extender fragment over $\cR'$
for some ordinal $\gm$ and premouse $\cM$ such that
$F$ is its own trivial completion and
satisfies the coherence condition:
$$\Edot^{ult_0 ( \cM , F)} \res (\ga + 1) = \Evec \res \ga$$
and the initial segment condition:
if $F$ has a last generator, $\eta$, and the trivial completion of 
$G$ of $F\res \eta$ has length $\gg$,
then either $G \subseteq \Edot^\cR_\gg$ or
$$G \subseteq \Edot^{ult_0 (J^\cR_\eta , \Edot^\cR_\eta )}_\gg \ .$$
\vskip 3ex

So a protomouse is a premouse if and only if its last predicate
is a total extender over its universe, and not merely an extender fragment.
The definitions and results in Sections 1-4 of \cite{MiSt}
apply as well to protomice as to premice.  The next definition
uses notions from \cite{MiSt} as applied to protomice.
As is the case with premice,
if $\cR$ is a type III protomouse, then whenever we write
$\cR$, we almost always mean $\cR^{sq}$.

\vskip 3ex
\noindent
{\sc Definition 2.3.3.}
Suppose that $\cR$ is a protomouse and $\gk$ is a cardinal in $\cR$.
If there is a largest integer $n$ such
that $\gr^\cR_n > \gk$,
then we put $n(\cR , \gk) = n$;
otherwise, set $n(\cR , \gk ) = \go$.
Say $n(\cR , \gk ) = n < \go$.
Suppose that $\cR$ is $n$-sound;
let $u = u_n (\cR )$ = the solidity witness for $p_n (\cR )$
and $p = p_{n+1}^{(\cR , u )}$ =
the $(n+1)$'st standard parameter of $(\cR , u)$.
Suppose that $p - \gk = \{ p_0 > \cdots > p_k \}$ and 
whenever $i \leq k$, the $i$'th solidity witness for $p$ is in $\cR$,
that is,
$$w_i = Th_{n+1} ( p_i \ \cup \ p|i ) \in | \cR |.$$
Then we say that $\cR$ is {\bf solid above $\gk$}
and we set 
$p (\cR , \gk) = \langle p,u, \langle w_i | i \leq k \rangle \rangle$.
If, in addition,
$$| \cR | = H^\cR_{n+1} ( \gk \ \cup \  p(\cR ,\ga)  ),$$
then we say that $\cR$ is {\bf $\gk$-sound}.

\vskip 6ex
\noindent
{\sc {\S}2.4 \ \underline{Iteration Trees on Phalanxes of Protomice}}
\vskip 4ex

Our definition of a phalanx will differ from both definition 6.5 of
\cite{St1} and definition 9.6 of \cite{St1}, although the flavor is very
much the same.
\vskip 3ex

\noindent
{\sc Definition 2.4.1.}
Suppose that $\glvec = \langle \gl_\ga \, | \, \ga < \gg \rangle$
is an increasing sequence of ordinals,
and $\cRvec = \langle \cR_\ga \, | \, \ga \leq \gg \rangle$
is a sequence of protomice.
Then $(\cRvec , \glvec )$ is a 
{\bf phalanx of protomice of length $\gg$}
if and only if whenever $\ga \leq \gb \leq \gg$,
$\cR_\ga$ and $\cR_\gb$ agree below $\gl_\ga$.
\vskip 3ex

To us, a phalanx will always be a phalanx of protomice.
Note that in Definition 2.4.1, 
we do not require that $\cR_\ga$ be a premouse,
nor that $\gl_\ga$ be a cardinal in $\cR_\gb$.
$\cR_\gg$ is called the {\bf starting model}.
We often write $((\cRvec \res \gg , \cR_\gg ), \glvec )$
for $(\cRvec , \glvec )$.
Or, if $\gg = 2$, we write $((\cR_0 , \cR_1 ), \gl_0 )$.

\vskip 3ex
\noindent
{\sc Definition 2.4.2.}
Let $G$ be anything
and suppose that $\cR$ is a protomouse such that 
$$F = G \cap ( [ \gO ]^{< \go} \times | \cR | )$$
is an extender over $\cR$ with $\crit (F) = \gk$.
Under these conditions,
we define $\ult (\cR , G)$
to be $\ult_n (\cR , F)$,
where
$n = n (\cR , \gk )$.
\vskip 3ex

Iteration trees on structures similar to phalanxes of protomice 
were defined in \cite{MiSt} (where they were called ``pseudo iteration trees''),
and in Sections 6 and 9 of \cite{St1}.
Using Definition 2.4.2,
we can modify these earlier definitions in the obvious way to allow a phalanx
of protomice as an allowable base for an iteration tree.

\vskip 3ex
\noindent
{\sc Definition 2.4.3.}
By an {\bf iteration tree on a phalanx $(\cR , \glvec )$}
we mean a padded, $\go$-maximal pseudo-iteration tree
on the phalanx of protomice $(\cR , \glvec )$.
\vskip 3ex

There are several implicit assumptions in Definition 2.4.3
which we clear up here.
Suppose that $(\cR , \glvec )$ is a phalanx of length $\gg$
on which $\cT$ is an iteration tree.
Suppose that $\gg + \xi$ is least such that 
$E^\cT_{\gg + \xi} \not= \emptyset$.
It is conceivable that $\gn^\cT_{\gg +\xi} = \gn ( E^\cT_{\gg + \xi} )
< \gl_\ga$ for some $\ga < \gg$;
suppose that this is the case.
Let $\ga$ be minimal with this property.
If $\gg + \xi + 1 < \gg + \eta + 1 < \lh ( \cT )$ and
$\gn ( E^\cT_{\gg + \xi} )
< \crit ( E^\cT_{\gg + \eta } )
< \gl_\ga$, then it is not entirely clear from Definition 2.4.3
whether it is $\ga$ or $\gg$ 
that is the $\cT$-predecessor of $\gg + \eta + 1$.
The answer is that $\gg = \pred^\cT (\gg + \eta + 1 )$.
The relevant rule is that in applying 
$E^\cT_{\gg + \xi}$, we go back only as far as needed
to avoid moving generators, and only then consider
models before the starting model on the phalanx.
We remark that in practise, $\gn^\cT_{\gg + \xi}$
will be $\geq \gl_\ga$ for every $\ga < \gg$.

Another implicit assumption in Definition 2.4.2 is
that $(\cM^*_{\gg + \eta + 1})^\cT$ is\\
$\deg^\cT (\gg + \eta + 1)$-sound
whenever $E^\cT_{\gg + \xi} \not= \emptyset$.
In practise, $(\cM^*_{\gg + \eta + 1})^\cT$ will also be
solid above $\crit (E^\cT_{\gg + \eta + 1})$.

\vskip 3ex
\noindent
{\sc Definition 2.4.4.}
If $(\cRvec , \glvec)$ is a phalanx of length $\gg$
on which $\cT$ is an iteration tree,
and $\eta < \lh ( \cT )$, then by {\bf root$^\cT (\eta)$}
we mean the unique $\ga \leq \gg$ such that
$\ga \leq_T \eta$.
\vskip 3ex

The next definition is made relative to our choice of $\gU$.
The notation in the last two sections 
should be recalled.  A {\bf set protomouse}
is a protomouse of height $< \gO$; of course, every protomouse of height
$\gO$ is passive, and hence a weasel.

\vskip 3ex
\noindent
{\sc Definition 2.4.5.}
A phalanx $(\cRvec , \glvec)$ of length $\gg$ is {\bf special}
if and only if
there is a sequence
$\gkvec = \langle \ \gk_\ga \ | \ \ga < \gg \ \rangle$
such that:
\begin{list}{}{}
\item[({\it i})]{
	If $\ga \leq \gb < \gg$, then 
	$\gk_\ga$ is a cardinal in $\cR_\gb$.}
\item[({\it ii})]{
	If $\ga < \gg$, then
	$\gl_\ga = (\gk_\ga^+)^{\cR_\ga}$}
\item[({\it iii})]{
	If $\ga < \gg$, and $\cR_\ga$ is a set protomouse,
	then $\cR_\ga$ is $\gk_\ga$-sound.}
\item[({\it iv})]{
	If $\ga < \gg$, and $\cR_\ga$ is a weasel,
	then $\gk ( \cR_\ga ) \leq \gk_\ga$.}
\item[({\it v})]{
Suppose that
$\gb \leq \gg$ and $\cR_\gb$ is not a premouse.
Then there is a unique $\ga < \gb$
$\Fdot^{\cR_\gb}$
is a
$(\gk_\ga , \cR_\ga , \OR^{\cR_\gb})$-extender fragment over
$\cR_\gb$.\\
(That is, a unique $\ga < \gb$ such that $\Fdot^{\cR_\gb}$
is a
$(\gk_\ga , \OR^{\cR_\gb})$-extender over
$\cR_\ga$.)}
\end{list}

\vskip 3ex
\noindent
{\sc Definition 2.4.6.}
We say that $\cT$ is a {\bf special} iteration tree
if and only if 
there is a phalanx $(\cRvec , \glvec )$ of length $\gg$
such that $\cT$ is an iteration tree on $(\cRvec , \glvec )$,
and:
\begin{itemize}
\item	$(\cRvec , \glvec )$ is special.
\item	If $\gg + \xi$ is least such that 
	$E^\cT_{\gg + \xi} \not= \emptyset$,
	and $\ga < \gg$,\\
	then
	$\gk_\ga$ is a cardinal in 
	$J^{\cR_\gg}_{lh (E^\cT_{\gg + \xi } )}$.
\item	If $\ga = \pred^\cT ( \gg + \eta + 1) < \gg$,
	then $\crit (E^\cT_{\gg + \eta}) = \gk_\ga$.
\end{itemize}
\vskip 3ex

Suppose that 
$(\cRvec , \glvec)$ is a special phalanx of length $\gg$
on which $\cT$ is a special iteration tree.
Suppose that $\gb < \gg$ and 
$\cR_\gb$ is not a premouse.
By Definition 2.4.5({\it v}), there is an ordinal $\ga < \gb$ such that
$\Fdot^{\cR_\gb}$ is a 
$(\gk_\ga , \OR^{\cR_\gb})$-extender over
$\cR_\ga$.
Suppose that $\gb = \rt^\cT ( \gg + \eta )$ and that
$\cD^\cT \cap \, (\gb , \gg + \eta ]_T = \emptyset$.
Since $\gk_\ga < \gk_\gb = \crit (i^\cT_{\gb , \gg + \eta})$,
the iteration of $\cR_\gb$ out to $\cM^\cT_{\gg + \eta}$ stretches the
last extender fragment, without changing the critical
point or the sets that it measures.
So $\Fdot^{\cM^\cT_{\gg +\eta}}$ is also an
extender over
$\cR_\ga$ with critical point $\gk_\ga$, and 
$\Fdot^{\cM^\cT_{\gg +\eta}}$ is
a reasonable candidate for $E^\cT_{\gg + \eta}$
with $\ga = \pred^\cT ( \gg + \eta + 1)$.

Continuing the same discussion, 
it is also possible that $E^\cT_{\gg +\eta}$ is an extender over an
initial segment of $\cM^\cT_{\gg + \eta}$,
but that it measures more subsets of its critical point than are in
the model to which it is being applied.
This happens, for example, when  $\gb = \rt^\cT (\gg + \eta )$,
$E^\cT_{\gg + \eta}$ is an extender from the sequence
$\Edot^{\cM^\cT_{\gg + \eta}}$
with $\crit (E^\cT_{\gg + \eta} ) = \gk_\ga$ for some $\ga < \gb$,
so that $\ga = \pred^\cT (\gg + \eta +1)$,
and 
$\gl_\ga < (\gk_\ga^+ )^{\cM^\cT_{\gg + \eta}}$.
Notice that in this case,
$$\cM^\cT_{\gg + \eta +1}\mbox{ and }
\ult ( J^{\cM^\cT_{\gg + \eta}}_{lh (E^\cT_{\gg + \eta} )}, E^\cT_{\gg + \eta} )
\mbox{ agree below }
(\gn ( E^\cT_{\gg + \eta} )^+ )^{\cM^\cT_{\gg + \eta +1}},$$
but it is possible that the successor of
$\gn ( E^\cT_{\gg + \eta} )$
as computed in $\cM^\cT_{\gg + \eta +1}$
is strictly less than
the successor of
$\gn ( E^\cT_{\gg + \eta} )$
as computed in
$$\ult (J^{\cM^\cT_{\gg + \eta}}_{lh (E^\cT_{\gg + \eta} )}, 
E^\cT_{\gg + \eta} ).$$
This makes it unclear how we can compare
two arbitrary special phalanxes
(in the sense of Section 7 of \cite{MiSt}).

\vskip 3ex
\noindent
{\sc Definition 2.4.7.}
Suppose that $(\cRvec , \glvec)$ is a special phalanx of length $\gg$
with the associated sequence $\gkvec$ as above.
$(\cRvec , \glvec)$ is {\bf very special} if and only if
for $\ga < \gg$:
\begin{itemize}
\item	If $\cR_\ga$ is a weasel, 
	then $c(\cR_\ga ) - \gk_\ga = \emptyset$.
\item	If $\cR_\ga$ is not a premouse,
	then $\cR_\ga$ is an active protomouse of type I or II
	and $n ( \cR_\ga , \gk_\ga ) = 0$.
\end{itemize}
\vskip 3ex

Most of the phalanxes
$(\cRvec , \glvec)$  that come up in our proof of Theorem 1.1
are special, and some of these are very special.

\vskip 6ex
\noindent
{\sc {\S}2.5 \ \underline{Long Extenders}}
\vskip 4ex

Throughout this section, $X$ is an elementary substructure
of $V_{\gO +1}$ (together with the measure on $\gO$ as a predicate)
and $\gp : N \lra V_{\gO +1}$ is the inverse of 
the transitive collapsing map for $X$.  Let $\gd = \crit (\gp )$
and $E_\pi$ be the
$(\gd , \gO)$-extender derived from $\pi$ (see Steel \cite{St2}).
So
$$E_\pi = \langle (E_\gp )_a \ | \ a \in [\gO ]^{< \go} \rangle,$$
where for any $a \in [ \gO ]^{<\go}$,
if $\gd_a$ is the least $\gdtil \geq \gd$ such that $a \subset \gdtil$, then
$$(E_\gp )_a = \{ x \subseteq [ \gd_a ]^{|a|} \ | \ a \in \gp (x) \}.$$
The ordinals $\gd_a$ we call the {\bf critical points} of $E_\pi$;
note that there may be many.

Suppose that $\cP$ is a premouse and that $\gk$ is a cardinal
in $\cP$ such that
$\cP (\gk ) \cap |\cP | \subset N$
and $\cP \models \gk^+$ exists.
Let 
$$F = E_\gp \ \cap \ ( [\gp (\gk ) ]^{< \go} \times | \cP | ).$$
In this situation, we wish to define 
$\ult (\cP , E_\gp \res \gp (\gk ) )$.

First suppose that $\cP$ is either passive or active of type I or II.
Let $n = n(\cP , \gk )$.
Then $\ult (\cP , E_\gp \res \gp (\gk ) )$
is $\ult_n (\cP , F)$, the ultrapower of $\cP$ by $F$
using certain functions $f : [ \gk ]^k \lra |\cP |$, where $k < \go$.
In $n=0$, we use only $f \in |\cP |$.
If $0 < n < \go$, 
then functions $f$ of the form $u \mapsto \gt^\cP_\varphi [ u , q]$
for a $\gS_n$-Skolem term $\gt_\varphi$ and parameter $q$ from $\cP$
are also allowed.
If $n=\go$, then all $f$ that
are first order definable over $\cP$ are allowed.

Now suppose that $\cP$ is active of type III. 
In this case, whenever we make reference to the fine structure of 
$\cP$, we are really referring to the fine structure of $\cP^{sq}$.
So, for example, $\gr_n^\cP$ is really $\gr_n^{\cP^{sq}}$.
With this interpretation, let
$n$ be as above.
It is not difficult to extend the results in Section 3 of \cite{MiSt}
to see that $\ult_n (\cP^{sq} , F)$ is itself a squashed potential premouse
(sppm).
By 
$\ult (\cP , E_\gp \res \gp (\gk ) )$,
we mean the premouse $\cR$ such that
$\cR^{sq} = \ult_n (\cP^{sq} , F)$.

In either case, let $\cR =
\ult (\cP , E_\gp \res \gp (\gk ) )$.
We assume throughout this section that $\cR$ is well-founded.
Let $\gptil$ be the ultrapower map.  So,
except when $\cP$ is active of type III, 
the map $\gptil : \cP \lra \cR$ is given by 
$$\gptil ( x ) = \left[ a,u\mapsto x \right]^\cP_{E_\pi \res \gp (\gk )}.$$
And when $\cP$ is active of type III,
the map $\gptil : \cP^{sq} \lra \cR^{sq}$ is given by 
$$\gptil ( x ) = \left[ a,u\mapsto x 
\right]^{\cP^{sq}}_{E_\pi \res \gp (\gk )} \ .$$

\vskip 3ex
\noindent
{\sc Lemma 2.5.1.}
Suppose that $\cP$ above is a $\gk$-sound premouse 
and that $n = n(\cP , \gk ) < \go$.
Suppose further that if $n=0$, $\cP$ is active of type I or II,
and $\gmdot^\cP < \gk$,
then
$$\sup \left( \gp '' \left( \gmdot^+ \right)^\cP \right)
= \gp \left( \left( \gmdot^+ \right)^\cP \right).$$
Then
\begin{list}{}{}
\item[1.]{$\gptil \res ( \gk^+ )^\cP = \gp \res (\gk^+ )^\cP$.}
\item[2.]{$\cR$ is a premouse of the same type as $\cP$.} 
\item[3.]{
$\gptil$ is an $n$-embedding that is continuous at
$\gl$ whenever $\gk < \gl < \gr_n^\cP$ and $\cP \models \cf (\gl ) > \gk$.}
\item[4.]{$\cR$ is $\gp (\gk )$-sound.}
\item[5.]{$\gptil (p(\cP , \gk )) = p(\cR , \gp (\gk ))$.}
\end{list}
\vskip 3ex

\noindent
{\sc Lemma 2.5.2.}
Suppose that $\cP$ above is a $\gk$-sound premouse that is active of type I 
or II with 
$n(\cP , \gk ) = 0$, $\gmdot^\cP < \gk$,
and 
$$\sup \left( \gp '' \left( \gmdot^+ \right)^\cP \right)
< \gp \left( \left( \gmdot^+ \right)^\cP \right).$$
Then
\begin{list}{}{}
\item[1.]{$\gptil \res ( \gk^+ )^\cP = \gp \res (\gk^+ )^\cP$.}
\item[2.]{$\cR$ is {\it not} a premouse;  it is a protomouse of the same
type as $\cP$.}
\item[3.]{
$\gptil$ is $\gS_1$-elementary,
$\sup ( \gptil '' \OR^\cP ) = \OR^\cR$,
and $\gptil$ is continuous at
$\gl$ whenever $\gk < \gl < \OR^\cP$ and $\cP \models \cf (\gl ) > \gk$.}
\item[4.]{$\cR$ is $\gp (\gk )$-sound.} 
\item[5.]{$\gptil ( p (\cP , \gk ) ) = p (\cR , \gp ( \gk ))$.}
\end{list}
\vskip 3ex

\noindent
{\sc Lemma 2.5.3.}
Suppose that $\cP$ above is a weasel with $\gk(\cP ) \leq \gk$.
Then
\begin{list}{}{}
\item[1.]{$\gptil \res ( \gk^+ )^\cP = \gp \res (\gk^+ )^\cP$.}
\item[2.]{$\cR$ is a weasel.}
\item[3.]{
$\gptil$ is fully elementary and is continuous
at $\gl$ whenever $\gk < \gl < \gO$ and $\cP \models \cf (\gl ) > \gk$.}
\item[4.]{$\gk ( \cR ) \leq \gp (\gk )$}
\item[5.]{$c(\cR ) = \gptil (c ( \cP ))$.}
\end{list}
\vskip 3ex

\noindent
{\sc Lemma 2.5.4.}
Suppose that $\cP$ is a set premouse with
$\gtdot^\cP \leq \gk$ and that $\cP$ is Dodd-solid above $\gk$.
Then $\gtdot^\cR \leq \pi ( \gk )$, $\cR$ is Dodd-solid above $\pi (\gk )$, 
and
$$\{ \sdot^\cR_i \ | \ \sdot^\cR_i \geq \pi ( \gk ) \} = 
\gptil \left( \{ \sdot^\cP_i \ | \ \sdot^\cP_i \geq \gk \} \right) \ .$$
\vskip 3ex

\noindent
{\sc Lemma 2.5.5.}
Suppose that $\cP'$ is a protomouse
such that $\cP$ and $\cP'$ agree below
$( \gk^+ )^\cP$, 
and 
$\cR' = \ult ( \cP' , E_\gp \res \gp ( \gk ))$ 
is well-founded.
Then $\cR'$ and $\cR$ agree below
$$(\pi (\gk)^+)^\cR
= \pi( (\gk^+)^\cP )
= \sup (\pi '' (\gk^+)^\cP )$$
$$= \{ \gp (f)(a) \ 
| \ a\in [\gp( \gk ) ]^{<\go} \ \& \ f: [\gk ]^{|a|} \ \ra ( \gk^+ )^\cP
\ \& \ f \in | \cP | \} .$$

\vskip 6ex
\noindent
{\sc {\S}2.6 \ \underline{Lifting by a Dodd-Squashed Extender Fragment}}
\vskip 4ex

The proof of the next lemma is implicit in the proof of weak square
in Section 5 of \cite{Sch}.

\vskip 3ex
\noindent
{\sc Lemma 2.6.1.}
Let $\cR$ be a protomouse, $\gm = \gmdot^\cR$, $F^* = \Fdot^\cR$, 
$\ga = \lh (F) = \OR^\cR$, $s = \sdot^\cR$, and $\gt = \gtdot^\cR$.
Let $\cM$ be a protomouse such that
$F$ is a $(\gm ,\ga)$-extender over $\cM$.
Suppose that $\gk$ is a cardinal in $\cR$ such that $\gt \leq \gk$
and $\cR$ is Dodd-solid above $\gk$.
Let $\cS = \ult (\cM , F)$ and let $j: \cM \lra \cS$ be the ultrapower
embedding.
So $\cS = \ult (\cM , F \res (\gk \cup s))$,
$\decap ( \cR )$ is a proper initial segment of $\cS$,
and $\OR^\cR$ is the cardinal successor of $\gndot^\cR$ in $\cS$.

Given the above, the following hold.

\begin{list}{}{}
\item[A.]{Suppose that $\cM$ is $\gm$-sound
with $p (\cM , \gm ) = \langle p , u , w \rangle$.
Then $\cS$ is $\gk$-sound
with $p(\cS , \gk ) = \langle j(p) \cup s , j(u) , w' \rangle$
for the right $w'$.}
\vskip 2ex
\item[B.]{Suppose that $\cM$ is a weasel
with $\gk (\cM) \leq \gm$.
Then $\cS$ is a weasel with
$\gk(\cS ) \leq \gk$ and $c(\cS ) - \gk = j(c(\cM )) \cup s$.
In particular, 
$\cS$ has the $c(\cS)$-hull property and the $c(\cS)$-definability
property above $\gk$.}
\end{list}

\vskip 6ex
\noindent
{\sc {\S}3 \ \underline{Proof of the Weak Covering Property}}
\vskip 4ex

Assume the hypothesis of Theorem 1.1.
So $\gk < \gO$ is a cardinal of $K$, 
$\card (\gk )$ is countably closed,
$\gO$ is measurable, and there is no inner model with a Woodin cardinal.
Construct $K$ and $K^c$ as subsets of $V_\gO$ using a fixed measure on 
$\cP (\gO )$.
Let $\gl = (\gk^+ )^K$ and let $\gU$ be any cardinal $> \gl$.
Let $X$ be an elementary substructure of $V_{\gO +1}$ 
(together with a predicate for the measure on $\gO$)
such that:
$$^\go (X \cap V_\gO ) \subseteq X$$
$$\card (X) < \card ( \gk )$$
$$\gk , \gl , \gU , \gO \in X  \ .$$
Such $X$ exist because $\card (\gk )$ is countably closed.
Note that if $\cf (\gl ) < \card (\gl )$, then there would be such $X$
with the additional property that 
$\gl \, \cap X$ is cofinal in $\gl$; in fact, this could be achieved with 
$$\card (X) = ( \cf(\gl ) )^\go \cdot \go_1 < \card ( \gk ) \ .$$
But we prefer to argue directly, rather than by contradiction,
so we do not assume, at least for now, that $\cf (\gl ) < \card (\gl )$.
We shall develop properties of $X$ (which are of independent interest),
and then somewhat later
argue by contradiction that Theorem 1.1 holds (just after Corollary 3.12).

Let $\gp : N \lra V_{\gO +1}$ be the inverse of the transitive collapse
of $X$ and let $\gd = \crit (\gp )$.
Let $W$ be the canonical very soundness witness for $\cJ^K_{\gU   }$.
and let $\Wbar = \gp^{-1} (W)$.
Compare
$W$ and $\Wbar$ as in Section 7 of \cite{MiSt}.  The result is an ordinal
$\theta \leq \gO$ and $\go$-maximal, padded iteration trees $\cT$ on $W$ and 
$\cTbar$ on $\Wbar$, both of length $\theta +1$.  Note
that $\Wbar$ is $\gO$-iterable, as any putative iteration tree
on $\Wbar$ can be copied to a putative iteration tree on $W$ via $\gp$,
and, of course, $W$ is $\gO$-iterable.
For $\eta \leq \theta$, put 
$W_\eta = \cM^\cT_\eta$ and $\Wbar_\eta = \cM^\cTbar_\eta$.
Because $W$ is universal, $\Wbar_\theta$ is an initial segment of 
$W_\theta$ and there is no dropping along the main branch 
of $\cTbar$, i.e.,
$$[0,\theta ]_{\Tbar} \ \cap \ \cD^\cTbar = \emptyset.$$
We shall see, in fact, that $\Wbar = \Wbar_\theta$.

Let $\gkbar = \gp^{-1} ( \gk )$, $\glbar = \gp^{-1} (\gl )$,
$\gUbar = \gp^{-1} ( \gU )$, and $\gObar = \gp^{-1} ( \gO )$.
For $\eta \leq \theta$, put
$E_\eta = E^\cT_\eta$ and $\Ebar_\eta = E^\cTbar_\eta$.
Choose $\gkvec$ and $\gg$ so that
$$\gkvec = \langle \ha^{\Wbar_\theta}_\ga \, | \, \ga \leq \gg \rangle$$
lists the cardinals of $\Wbar_\theta$ that are $\leq \gUbar$
and let 
$$\glvec = \langle \gk_{\ga +1} \, | \, \ga < \gg \rangle \ .$$
We define premice $\cP_\ga$ and ordinals $\eta (\ga )$ for 
$\ga \leq \gg$ as follows
If there is an $\eta$ such that $\gk_\ga < \nu^\cT_\eta$,
then let $\eta (\ga )$ be the least such $\eta$;
otherwise, let $\eta (\ga ) = \theta$.
If there is an initial segment $\cP$ of $W_{\eta ( \ga )}$
with $\OR^\cP \geq \gl_\ga$ and $\gr_\go^\cP \leq \gk_\ga$,
then let $\cP_\ga$ be the shortest such $\cP$;
otherwise, let $\cP_\ga = W_{\eta (\ga )}$.
Observe that $\cP_\ga$ is exactly the premouse to which
an extender with critical point $\gk_\ga$ would be applied in
an iteration tree extending $\cT$.
When $\gl_\ga < \gd$, $\cP_\ga$ is just $W$.
Let $\cPvec = \langle \cP_\ga \, | \, \ga \leq \gg \rangle$.
Then $(\cPvec , \glvec )$ is a very special phalanx of premice.
Moreover, any iteration tree on $(\cPvec , \glvec )$
can be construed as an iteration tree extending $\cT$,
so $(\cPvec , \glvec )$ is iterable.

For $\ga \leq \gg$, let 
$\cR_\ga = \ult (\cP_\ga , E_\gp \res \gp (\gk_\ga ))$,
as long as this ultrapower is wellfounded.
Also, let $\gp_\ga : \cP_\ga \lra \cR_\ga$ be the
ultrapower map and $m_\ga = n (\cP_\ga , \gk_\ga)$.
For $\ga < \gg$, put
$\gL_\ga = \sup ( \gp_\ga \, '' \gl_\ga ) 
= (\gp_\ga (\gk_\ga )^+ )^{\cR_\ga}$.
So $\cRvec = \langle \, \cR_\ga \, | \, \ga \leq \gg \, \rangle$
and $\gLvec = \langle \, \gL_\ga \, | \, \ga < \gg \, \rangle$.
Our first essential use of countable closure is the following lemma;
the other will be Lemmas 3.3 and  3.13.

\vskip 3ex
\noindent
{\sc Lemma 3.1.}
For every $\ga \leq \gg$, the ultrapower defining $\cR_\ga$
is wellfounded. Moreover,
$(\cRvec , \gLvec )$
is an iterable, very special phalanx.
\vskip 3ex

\noindent
{\sc Proof sketch.}
The well-foundedness and iterability follow from
the countable completeness of $E_\pi$ in the standard way.
The argument is very similar to the one we shall give
in the proof of Lemma 3.13.
The rest follows from
the lemmas in {\S}2.5.\\
\qed {\tiny \ \ Lemma 3.1}
\vskip 3ex

We give some motivation for what is to come.
The basic idea behind earlier proofs that 
$K$ satisfies the weak covering property lower down in the
large cardinal hierarchy (i.e., those in \cite{DoJe2}, \cite{Mi2},
and \cite{Je2})
is as follows.
First, show that
$\Wbar$ does not move in its comparison with $W$;
this would imply that $\gkbar = \gk_{\gg_0}$ for some $\gg_0 < \gg$,
and that $( \gk^+ )^{\cR_{\gg_0}} = \gl$
Second, show that $\cR_{\gg_0}$ is a set premouse;
then $\cR_{\gg_0}$ would be $\gk$-sound.
Third, show that $\cR_{\gg_0}$ is iterable in a way that implies
that $\cR_{\gg_0}$ is an element of $K$,
which would immediately give a contradiction.

We shall be able to carry out a version of the first step.
By induction on $\ga \leq \gg$, we shall show the following;
we remark that this is the first of several induction hypotheses to come.

\vskip 3ex
{\bf (1)$_\ga$} \ \ If $\Ebar_\eta \not= \emptyset$,
then $\lh ( \Ebar_\eta ) > \gl_\ga$.
\vskip 3ex

However, there does not seem to be much hope of carrying out the second
step, that is, in showing that $\cR_{\gg_0}$
is a set premouse.
But if $\cR_{\gg_0}$ is a weasel, then clearly it is not an element of $K$.
And if $\cR_{\gg_0}$ is merely a protomouse
that is not a premouse, then it is not clear what sort of iterability
of $\cR_{\gg_0}$ would imply that it is in $K$.
So our next step will be to identify a phalanx
of premice
$(\cSvec , \gLvec )$ with properties very similar to those of
$(\cRvec , \gLvec )$.
Later, in order to establish some connection between $\cS_{\gg_0}$ and $K$,
we shall show inductively that $\cS_\ga$ is iterable in a strong sense.
In the end, we do not show that $\cS_{\gg_0}$
is a set premouse, but rather just deal with the possibility that it is
a weasel.

We define $\cS_\gb$ by induction on $\gb \leq \gg$ as follows.
If $\cR_\gb$ is a premouse, then $\cS_\gb = \cR_\gb$.
Suppose that $\cR_\gb$ is not a premouse.
Then, since $(\cRvec , \gLvec )$ is a special phalanx,
there is an ordinal
$\ga < \gb$ such that 
$\pi ( \gk_\ga ) = \gmdot^{\cR_\gb} = \crit (\Fdot^{\cR_\gb})$.
Let 
$\cS_\gb = \ult (\cS_\ga , \Fdot^{\cR_\gb } )$,
so long as this ultrapower is wellfounded.
Before showing that this ultrapower is always wellfounded,
we define some related premice $\cQ_\gb$ by induction
on $\gb \leq \gg$.  If $\cS_\gb = \cR_\gb$, then 
$\cQ_\gb = \cP_\gb$.  Suppose $\cS_\gb \not= \cR_\gb$.
Then there is an ordinal $\ga < \gb$ such that 
$\gk_\ga = \gmdot^{\cP_\gb}$.
Let $\cQ_\gb = \ult (\cQ_\ga , \Fdot^{\cP_\gb} )$,
so long as this ultrapower is wellfounded.
The proof of the following fact is 
is similar to the proof of the main result in Section 9 of \cite{St1};
we omit the details.

\vskip 3ex
\noindent
{\sc Fact 3.2.}
For every $\gb \leq \gg$, the ultrapower defining $\cQ_\gb$
is wellfounded.

\vskip 3ex
\noindent
{\sc Lemma 3.3.}
For every $\gb \leq \gg$, the ultrapower defining $\cS_\gb$
is wellfounded.
\vskip 3ex

\noindent
{\sc Proof sketch.}
The proof of Lemma 3.3 is essentially included in the proof
of Lemma 3.13 to come.  Very briefly, the idea is as follows.
First show
that $\cS_\gb$ is the ultrapower
of $\cQ_\gb$ by $E_\pi \res \pi ( \gk_\gb )$.
And then use the countable completeness of $E_\pi$
in the usual way to see that $\cS_\gb$ is wellfounded.
\qed {\tiny \ \ Lemma 3.3}
\vskip 3ex

We want to show that $(\cSvec , \gLvec )$ is a special phalanx of premice.
But first we need some preliminary lemmas.

\vskip 3ex
\noindent
{\sc Lemma 3.4.}
Suppose that $\cR_\ga$ is not a premouse.
Then $\gtdot^{\cP_\ga} \leq \gk_\ga$
and $\cP_\ga$ is Dodd-solid above $\gk_\ga$.
\vskip 3ex

\noindent
{\sc Proof.}
Since $\cR_\ga$ is not a premouse,
the lemmas in {\S}2.5 imply that $\cP_\ga$ is a $\gk_\ga$-sound
premouse that is active of type I or II,
$n(\cP_\ga , \gk_\ga ) = 0$,
$\gmdot^{\cP_\ga} < \gk_\ga$ and $\gp$ is discontinuous at 
$(\gmdot^+ )^{\cP_\ga}$.
If $\cP_\ga$ is type I, then 
$\gk_\ga = (\gmdot^+ )^{\cP_\ga} = \gtdot^{\cP_\ga}$
and $\sdot^{\cP_\ga} = \emptyset$, so we are done.
Suppose that $\cP_\ga$ is type II, but that the conclusion
of Lemma 3.4 fails.
Then, by Corollary 2.1.6,
$(\gmdot^+ )^{\cP_\ga}$ has countable cofinality.
But $^{\go} (X \cap V_\gO ) \subset X$, so $\gp$ is continuous at 
$(\gmdot^+ )^{\cP_\ga}$.  Contradiction.
\qed {\tiny \ \ Lemma 3.4}
\vskip 3ex

We do not consider Lemma 3.4 and essential use of countable closure,
since the proof only uses the continuity of $\gp$ at ordinals of countable
cofinality, which is easy to arrange.
By combining Lemmas 3.4 and 2.5.4, we get the following corollary.

\vskip 3ex
\noindent
{\sc Corollary 3.5.}
If $\cR_\ga$ is not a premouse,
then $\gtdot^{\cR_\ga} \leq \gp ( \gk_\ga )$ and
$\cR_\ga$ is Dodd-solid above $\gp (\gk_\ga )$.
\vskip 3ex

Recall that $m_\ga = n ( \cP_\ga , \gk_\ga ) = n(\cR , \gp ( \gk_\ga ))$.
Define $n_\ga = 
n(\cQ_\ga , \gk_\ga )$.  Then also $n_\ga
= n(\cS_\ga , \gp ( \gk_\ga ))$.
The next lemma follows by an easy induction on $\ga \leq \gg$,
using Lemma 2.6.1.  The analogous statement (really part
of the proof) holds 
for the $\cP_\ga$'s, $\cQ_\ga$'s, and $\gk_\ga$'s.

\vskip 3ex
\noindent
{\sc Lemma 3.6.}
Consider any $\ga \leq \gg$.
There is a unique parameter
$\{ \ga_0 > \ga_1 > \cdots > \ga_k \}$ 
with
\begin{list}{}{}
\item[1.]{$\ga_0 = \ga$,}
\item[2.]{$\gp ( \gk_{\ga_i} ) = \gmdot^{\cR_{\ga_i}}  $ whenever $0 \leq i
\leq k$, and}
\item[3.]{$\cR_{\ga_k}$ is a premouse.}
\end{list}
We have the sequence of ultrapower maps:
$$\cR_{\ga_k} = \cS_{\ga_k} \lra \cS_{\ga_{k-1}} \lra \cdots \lra
\cS_{\ga_0} = \cS_\ga$$ 
which we call the {\bf decomposition} of $\cS_\ga$.
Let $\gs_i : \cS_{\ga_i} \lra \cS_\ga$ whenever $0 \leq i \leq k$,
with $\gs_0 = \id$.
Then the following hold.
\begin{list}{}{}
\item[(a)]{Suppose that $\cR_{\ga_k}$ is a weasel.
Then $\cS_\ga$ is a weasel with
$\gk (\cS_\ga ) \leq \pi ( \gk_\ga )$
and $$c(\cS_\ga ) - \pi ( \gk_\ga ) =
\bigcup_{0 < i \leq k} \gs_i ( \sdot^{\cR_{\ga_i}} ) \ .$$
In particular, $n_\ga = m_{\ga_k} = \go$.}
\item[(b)]{Suppose that $\cR_{\ga_k}$ is a set premouse
with 
$$p(\cR_{\ga_k} , \pi ( \gk_{\ga_k} ) ) =
\langle p,u,w \rangle \ .$$
Then $\cS_\ga$ is a $\pi ( \gk_\ga )$-sound set premouse
with $n_\ga = m_{\ga_k}$,
and for the right $u'$ and $w'$,
$$p(\cS_\ga , \pi ( \gk_\ga ) ) =
\langle p' ,u' ,w' \rangle$$
with 
$$p' = \gs_k (p) \cup
\bigcup_{0 < i \leq k} \gs_i ( \sdot^{\cR_{\ga_i}} ) \ .$$}
\end{list}

\vskip 3ex
\noindent
{\sc Definition 3.6.1.}
If $\cM$ is a protomouse,
then 
{\bf decap}$(\cM ) = \langle \, |\cM | \, , \, \in \, , \, \Edot^\cM \, 
\rangle$.
If $\xi \leq \OR^\cM$,
then $\cM \| \xi = \decap ( \cJ^\cM_\xi )$.
\vskip 3ex

An easy induction 
using, among other things, the coherence condition on $\Fdot^{\cR_\ga}$
gives the following lemma.  Again, the analogous statement is true of 
the $\cP_\ga$'s, $\cQ_\ga$'s and $\gk_\ga$'s.

\vskip 3ex
\noindent
{\sc Lemma 3.7.}
The notation here is as in Lemma 3.6.
Suppose that $\cR_\ga \not= \cS_\ga$.
Then:
\begin{itemize}
\item	$n_\ga = n_{\ga_0} = \cdots = n_{\ga_k} = m_{\ga_k} \leq \go$
\item	$m_\ga = m_{\ga_0} = \cdots = m_{\ga_{k-1}} = 0$
\item	$\OR^{\cR_\ga}$ is a successor cardinal in $\cS_\ga$
	and
	$\decap ( \cR_\ga )$
	is a proper initial segment of $\cS_\ga$.
\item	For $i = 0 , \dots , k-1$,
	$$\OR^{\cS_\ga}
	= \sup ( \gs_k \, '' \OR^{\cR_{\ga_k}} )
	\geq \gs_k ( \gr_{n_{\ga_k}} ( \cR_{\ga_k} ))
	= \gr_{n_\ga} ( \cS_\ga ) 
	> \gs_i ( \OR^{\cR_{\ga_i}} )
	\geq \OR^{\cR_\ga}$$
\end{itemize}

\vskip 3ex
\noindent
{\sc Corollary 3.8.}
Both $(\cSvec , \gLvec )$ and 
$(\cQvec , \glvec )$ are special phalanxes of premice.
\vskip 3ex

We remark that $(\cSvec , \gLvec )$ is an iterable phalanx;
this will follow from the proof of Theorem 1.1, 
but will not be isolated, nor
used. Next, we list induction hypotheses (2)$_\ga$ through (6)$_\ga$.
The reader should note the similarity between (2)$_\ga$ and
the definition of $\pi ( \gk_\ga )$-strong given in Section 6 of \cite{St1}.

\vskip 3ex
{\bf (2)$_\ga$} \ \ $((W, \cS_\ga ), \pi ( \gk_\ga ))$ is an iterable
phalanx of premice.

\vskip 2ex
{\bf (3)$_\ga$} \ \ $((\Wbar , \cQ_\ga ), \gk_\ga )$ is an iterable
phalanx of premice.

\vskip 2ex
{\bf (4)$_\ga$} \ \ $((\cPvec \res \ga , \Wbar ),\glvec \res \ga )$
is an iterable, very special phalanx of premice.

\vskip 2ex
{\bf (5)$_\ga$} \ \ $((\cRvec \res \ga , W ), \gLvec \res \ga )$ 
is a very special phalanx, that is iterable with respect to special iteration
trees.

\vskip 2ex
{\bf (6)$_\ga$} \ \ $(\cSvec \res \ga , W ), \gLvec \res \ga )$ 
is a very special phalanx, that is iterable with respect to special iteration
trees.
\vskip 3ex

Assuming that (1)$_\gb$ holds for
every $\gb < \ga$, we shall show:
$$\mbox{(3)$_\ga$ $\Lra$ (2)$_\ga$}$$
$$\mbox{(6)$_\ga$ $\Lra$ (5)$_\ga$ $\Lra$ (4)$_\ga$ $\Lra$ (1)$_\ga$}$$
$$\mbox{$\forall \gb < \ga$ (4)$_\gb$ $\Lra$ (3)$_\ga$}$$
$$\mbox{$\forall \gb < \ga$ (2)$_\gb$ $\Lra$ (6)$_\ga$}$$
Thus, by induction, (1)$_\ga$ through (6)$_\ga$ hold for all $\ga < \gg$.
What this buys us is explained by Corollary 3.12 below.

\vskip 3ex
\noindent
{\sc Lemma 3.9.}
Let $\ga \leq \gg$ be given.
Suppose that (2)$_\ga$ holds.
Then there is an iteration tree $\cU$ on $W$
with the following properties:
\begin{list}{}{}
\item[(a)]{$\cU$ has successor length, $\varphi + 1$.}
\item[(b)]{$\lh ( E^\cU_\eta )\geq \gp ( \gk_\ga )$ for every 
$\eta \leq \varphi$.}
\item[(c)]{
There is an initial segment $\cN$ of 
$\cM^\cU_\varphi$ and an $n_\ga$-embedding
$j: \cS_\ga \lra \cN$ such that 
$\crit (j) \geq \gp ( \gk_\ga )$.}
\end{list}

\vskip 3ex
\noindent
{\sc Lemma 3.10.}
Let $\ga \leq \gg$ be given.
Suppose that (1)$_\ga$ and 
(2)$_\ga$ hold and that $\cS_\ga$ is a set premouse.
If $E^W_{\gL_\ga} = \emptyset$, then $\cS_\ga$ is a proper
initial segment of $W$.
Otherwise, $\cS_\ga = \ult ( \cN , E^W_{\gL_\ga} )$
for the longest initial segment $\cN$ of $W$
over which $E^W_{\gL_\ga}$ is an extender.

\vskip 3ex
\noindent
{\sc Lemma 3.11.}
Let $\ga \leq \gg$ be given.
Suppose that (1)$_\ga$ and (2)$_\ga$ hold
and that $\cS_\ga$ is a weasel different from $W$.
Let $\cU$, $\varphi$, 
$j$, and $\cN$ 
be as in Lemma 3.9.  Then
\begin{list}{}{}
\item[(a)]{$j$ is an elementary embedding from $\cS_\ga$
into $\cN = \cM^\cU_\varphi$.}
\item[(b)]{$\lh ( E^\cU_\eta ) \geq \gL_\ga 
= ( \gp ( \gk_\ga )^+ )^{\cS_\ga}$ whenever 
$\eta \leq \varphi$.}
\item[(c)]{$1 \leq_U \varphi$}
\item[(d)]{The Dodd-projectum of $E^\cU_0$, $\gt (E^\cU_0 )$,
is $\leq \gp ( \gk_\ga )$;\\ 
in particular,
$\gL_\ga < (\gp ( \gk_\ga )^+ )^W$.}
\end{list}

\vskip 3ex
\noindent
{\sc Corollary 3.12.}
Suppose that $\ga \leq \gg$, $\gd \leq \gk_\ga$, 
and both (1)$_\ga$ and (2)$_\ga$ hold.
Then $\gL_\ga < ( \gp ( \gk_\ga )^+ )^W$.
Also, there is an iteration tree $\cU$ on $W$ satisfying 3.9(a) and 3.9(b)
and an elementary embedding $j : \cS_\ga \lra \cM^\cU_\varphi$
with $\crit (j) \geq \gp ( \gk_\ga )$
(if $\cS^\ga$ is a set, then $j$ is the identity).
\vskip 3ex

Theorem 1.1 will follow from Corollary 3.12
once we show that (1)$_\ga$ -- (6)$_\ga$ hold for all $\ga < \gg$.
Here is the argument.
Suppose that (1)$_\ga$ -- (6)$_\ga$ hold for all $\ga < \gg$.
Assume, for contradiction, that $\cf ( \gl ) < \card (\gk )$.
Then we may assume without loss of generality that $\gp$ is continuous at
$\glbar$.  
Since (1)$_\ga$ holds for all $\ga < \gg$,
there is some $\ga < \gg$ such that
$\gkbar = \gk_\ga$ and $\glbar = \gl_\ga$.
By the lemmas in \S2.5, $\gp \res \glbar = \gp_\ga \res \gl_\ga$.
Hence,
$$\gl = \gp ( \glbar ) = \sup ( \gp {\ }'' \glbar )
= \sup ( \gp_\ga  {\ }'' \gl_\ga )
= \gL_\ga \ .$$
Therefore
$\gl < (\gk^+ )^W = (\gk^+ )^K$,
by Corollary 3.12.
But this is a contradiction, since $\gl = (\gk^+ )^K$ by definition.

Just to be explicit, 
through the remainder of the proof,
we are not assuming that (1)$_\ga$ -- (6)$_\ga$
hold, but rather proving them
by induction on $\ga < \gg$.

Here are some further corollaries to Lemmas 3.9 and 3.10.
If (1)$_\ga$ and (2)$_\ga$ hold, $\cS_\ga$ is a set premouse,
and $\gp ( \gk_\ga )$ is a cardinal
in $W$, then
$\cS_\ga$ is a proper initial segment of $W$.
So if (2)$_\ga$ holds and (1)$_\gb$ holds for some $\gb \geq \ga$ such that
$\gk_\gb$ is a cardinal in $\Wbar$ (rather than just in $\Wbar_\theta$),
and $\cS_\ga$ is a set premouse, then $\cS_\ga$ is a proper initial segment
of $W$.

\vskip 3ex
\noindent
{\sc Proof of 3.9.}
Compare $W$ versus $((W , \cS_\ga ), \gp ( \gk ))$.
This can be done since we are assuming (2)$_\ga$.
This results in iteration trees $\cU$ on $W$ and
$\cV$ on 
$((W , \cS_\ga ), \gp ( \gk ))$.
For some $\varphi$,
$\lh (\cU ) = \varphi +1$ and $\lh (\cV ) = 1 + \varphi + 1$.
The comparison ends with either $\cM^\cU_\varphi$ 
an initial segment of 
$\cM^\cV_{1 + \varphi}$ or vice-versa.
Note that all extenders used on either iteration tree have 
length at least $\gp (\gk_\ga )$ since by (2)$_\ga$, 
$\cS_\ga$ and $W$ agree below $\gp (\gk_\ga )$.
The usual arguments using the universality,
hull, and definability properties of $W$ (as in Section 3 of \cite{St1})
show that $1 = \rt^V ( 1 + \varphi )$.
That is, the last model on $\cV$ lies above $\cS_\ga$ and not above
$W$.
Similarly, we see that $\cM^\cV_{1 + \varphi}$ is an initial segment of 
$\cM^\cU_\varphi$, that $[1 , 1+ \varphi ]_V \cap \cD^\cV = \emptyset$,
and that $\deg^\cV (1 + \varphi ) = n_\ga$.
Let $j = i^\cV_{1 , 1 + \varphi}$ and $\cN = \cM^\cV_{ 1 + \varphi }$ .
\qed {\tiny \ \ Lemma 3.9}
\vskip 3ex

\noindent
{\sc Proof of 3.10.}
We pick up where we left off in the proof of Lemma 3.9,
only now we assume that (1)$_\ga$ holds.
Then $\cS_\ga$ and $W$ agree below $\gL_\ga$,
so all extenders used on either of $\cU$ and $\cV$
have length at least $\gL_\ga$.
We also assume that $\cS_\ga$ is a set premouse.

\vskip 3ex
\noindent
{\sc Claim 1.}
$\cS_\ga$ does not move, that is, $\cN = \cS_\ga$.
\vskip 3ex

\noindent
{\sc Proof.}
Suppose otherwise.  Then $\cN$ is not 
$\gp (\gk_\ga )$-sound, so 
$[ 0 , \varphi ]_U \cap \cD^\cU \not= \emptyset$
and $\cM^\cU_\varphi = \cN$.
Let $\eta +1$ be the last drop (of any kind) along
$[0 , \varphi ]_U$. 
By this we mean that
$\eta + 1$ is largest
in $[0 , \varphi ]_U$ such that either $\eta +1 \in \cD^\cU$,
or $\deg^\cU (\eta +1 ) < \deg^\cU (\pred^\cU (\eta +1))$,
where we take
$\deg^\cU ( 0 ) = \deg^\cU ( 1 ) = \go$.
Notice that, by construction, $\xi +1$ is a drop whenever
$1 = \pred^\cU ( \xi +1 )$.
Let $\eta^* = \pred^\cU ( \eta +1 )$.

If $\crit (E^\cU_\eta) \geq \gp ( \gk_\ga )$,
then both $j$ 
and 
$i^\cU_{\eta  +1,\varphi}
\circ ( i^*_{\eta +1} )^\cU$
are the inverse of the transitive collapsing map
for 
$$H^\cN_{n_\ga +1} \left( \gp \left( \gk_\ga \right) \ \cup 
\ p\left( \cN , \gp \left( \gk_\ga \right)\right)\right) \ .$$
This gives the usual contradiction as in the proof of the
Comparison Lemma 7.1 of \cite{MiSt}.

So $\crit ( E^\cU_\eta ) < \gp ( \gk_\ga )$.
Therefore $\eta^* = 0$ and
$$n_\ga =
n(\cN , \gp ( \gk_\ga )) =
n(\cM^\cU_{\eta +1}, \gp ( \gk_\ga )) =
n( ( \cM^*_{\eta +1} )^\cU , \crit ( E^\cU_\eta )),$$
where $( \cM^*_{\eta +1} )^\cU$
is the proper initial segment of $W$ to which $E^\cU_\eta$ is applied.
Since 
$\cS_\ga$ and $\cN$ agree beyond $\gL_\ga = ( \pi ( \gk_\ga )^+)^{\cS_\ga}$
and 
$$\cS_\ga = 
\cH^\cN_{n_\ga +1} \left( \gp \left( \gk_\ga \right) \cup p \left( \cN , \gp \left( \gk_\ga \right) \right) \right) \ ,$$
we have that
$\cP ( \gp ( \gk_\ga )) \cap |\cN | \subseteq$
the transitive collapse of
$$H^\cN_{n_\ga +1} 
\left( \gp \left( \gk_\ga \right) \cup p \left( \cN , \gp \left( \gk_\ga
\right) \right) \right) \ .$$
But this is only possible if there are no generators
of $E^\cU_\eta$ $\geq 
\gp ( \gk_\ga )$.
So $\gn ( E^\cU_\eta ) = \gp ( \gk_\ga )$,
which in turn is only possible if $\eta = 0$ and 
$E^\cU_\eta = E^W_{\gL_\ga}$.
But then $\cM^\cU_{\eta +1}$ is $\gp ( \gk_\ga )$-sound,
so both $i^\cU_{\eta +1 , \varphi}$ and $j$
are the inverse of the transitive collapsing map for 
$$H^\cN_{n_\ga +1} \left( \gp \left( \gk_\ga \right) \ \cup 
\ p\left( \cN , \gp \left( \gk_\ga \right) \right) \right) \ .$$
This gives a contradiction as in the proof of the Comparison Lemma 7.1
of \cite{MiSt}.\\
\qed {\tiny \ \ Claim 1}
\vskip 3ex

Thus $\cS_\ga$ is an initial segment of $\cM^\cU_\varphi$.

\vskip 3ex
\noindent
{\sc Claim 2.}
If $\cS_\ga$ is a proper initial segment of $\cM^\cU_\varphi$,
then $\cM^\cU_\varphi = W$.
\vskip 3ex

\noindent
{\sc Proof.}
Suppose that $\cS_\ga$ is a proper initial segment of $\cM^\cU_\varphi$.
Since $\cS_\ga$ is $\gp ( \gk_\ga )$-sound,
$$\cS_\ga = \cJ^{\cM^\cU_\varphi}_\ge$$ for some 
$\ge < ( \gp ( \gk_\ga )^+ )^{\cM^\cU_\varphi}$.
But 
$( \gp ( \gk_\ga )^+ )^{\cM^\cU_\varphi}
\leq \lh (E^\cU_\eta )$ for the least $\eta$ 
such that $(\eta +1) \leq_U \varphi$.
So $\cS_\ga$ is an initial segment of $W$.
\qed {\tiny \ \ Claim 2}
\vskip 3ex

So we may assume that $\cS_\ga = \cM^\cU_\varphi$,
and therefore,
$[0 , \varphi ]_U \cap \cD^\cU \not= \emptyset$.
Since $\cS_\ga$ is $\gp ( \gk_\ga ) $-sound,
this is only possible if $\cS_\ga = 
\ult ((\cM_{\eta + 1}^* )^\cU , E^\cU_\eta )$,
for the least $\eta$ such that $(\eta +1) \leq_U \varphi$.
Reasons as before show that $E^\cU_\eta$ cannot have
generators $\geq \gp ( \gk_\ga )$.
Therefore $\eta = 0$ and $\lh ( E^\cU_\eta ) = \gL_\ga$.
\qed {\tiny \ \ Lemma 3.10}
\vskip 3ex

\noindent
{\sc Proof of 3.11.}
Again, we continue the comparison begun in the proof of Lemma 3.9.
Since we are assuming (1)$_\ga$,
$\cS_\ga$ and $W$ agree below $\gL_\ga$,
so all extenders used on either of $\cU$ and $\cV$
have length at least $\gL_\ga$.
We are also assuming also that $\cS_\ga$ is a weasel different from $W$.
So $\cS_\ga$ is a universal weasel with $\gk ( \cS ) \leq \gp ( \gk_\ga )$.
By the usual arguments using the universality, hull, and definability
properties of $W$ and $\cS_\ga$,
$\cM^\cU_\varphi = \cN$, 
$[0 , \varphi ]_U \cap \cD^\cU =
\emptyset$, and 
$[1 , 1 + \varphi ]_V \cap \cD^\cV =
\emptyset$, and
$j  : \cS_\ga \lra \cN$ is elementary.
Note that
$\cN$ has the $j (c ( \cS_\ga ))$-hull 
property at $\gp (\gk_\ga )$.
Let $\eta$ be least such that $(\eta +1) \leq_U \varphi$.

\vskip 3ex
\noindent
{\sc Claim 1.}
There is a parameter $t$ such that $\cM^\cU_{\eta +1}$
has the $t$-hull property at $\gp ( \gk_\ga )$.
\vskip 3ex

\noindent
{\sc Proof.}
Suppose that for every $t$,
the $t$-hull property fails at $\gp (\gk_\ga )$
in $\cM^\cU_{\eta +1}$.

\vskip 3ex
\noindent
{\sc Subclaim.}
For every
$\xi_0 \in [\eta +1 , \varphi ]_U$,
the $t$-hull property fails at $\gp (\gk_\ga )$
in $\cM^\cU_{\xi_0}$.
\vskip 3ex

\noindent
{\sc Proof.}
By induction on $\xi_0$.
Let $\xi_0$ be the least counterexample, and let $t$ be a witness.
Clearly $\xi_0 = \xi +1$ for some $\xi > \gg$.  Let 
$\xi^* = \pred^\cU ( \xi + 1 )$.  Say $t = i^\cU_{\xi^* , \xi + 1} ( f ) (a)$
with $a \in [ \gn ( E^\cU_\xi ) ]^{<\go}$,
$$f: [\crit (E^\cU_\xi ) ]^{|a|} \lra | \cM^\cU_{\xi^*} |,$$
and $f \in | \cM^\cU_{\xi^*} |$.
By Section 4 of \cite{MiSt}, $E^\cU_\xi$ is close to $\cM^\cU_{\xi^*}$.
In particular, $(E^\cU_\xi )_a $ is $\gS_1$ over $\cM^\cU_{\xi^*}$.
Since $\cM^\cU_{\xi^*}$ is a weasel,
$(E^\cU_\xi )_a \in |\cM_{\xi^*}^\cU |$.
By the induction hypothesis,
there is an $A \subseteq \gp ( \gk_\ga )$ such that 
$A \not\in$ the transitive collapse of 
$$H^{\cM^\cU_{\xi^*}}_\go \left( \gp ( \gk_\ga ) \cup 
\{ f , (E^\cU_\xi )_a \}  \cup \gG \right) \ ,$$
where $\gG$ = the fixed points of $i^\cU_{\xi^* , \xi + 1}$.
On the other hand, there is a term $\gt_\psi$ and $c\in [\gG ]^{< \go}$
such that 
$$A = \gt_\psi^{\cM^\cU_{\xi + 1}} [ t,c ] \cap \gp ( \gk_\ga ).$$ 
By {\L}o{\' s}' Theorem,
for $\gz < \gp ( \gk_\ga )$,
$$\gz \in A
\ \Llra \ \gz \in \gt_\psi^{\cM^\cU_{\xi^*}} [ f(u), c]
\mbox{ for } ( E^\cU_\xi )_a \mbox{ a.e. }u \ .$$
But then $A$ is not as it was chosen to be, a witness that the 
$\{ f , ( E^\cU_\xi )_a \}$-hull property fails 
at $\gp ( \gk_\ga )$ in $\cM^\cU_{\xi^*}$.
\qed {\tiny \ \ Subclaim}
\vskip 3ex

By the Subclaim, the $t$-hull property fails at $\gp ( \gk_\ga )$
in $\cN = \cM^\cU_\varphi$ for every $t$.
But this is a contradiction with $t = j (  c ( \cS_\ga ))$.
\qed {\tiny \ \ Claim 1}

\vskip 3ex
\noindent
{\sc Claim 2.}
$\crit (E^\cU_\eta ) < \gp (\gk_\ga )$.
\vskip 3ex

\noindent
{\sc Proof.}
We have two ``iterations'' from $W$ into $\cN$.
The first is $i^\cU_{0,\varphi}$.
The second is
$$W \stackrel{ i^\cT_{0,\eta (\ga_k )} }{\lra}
W_{\eta ( \ga_k )} = \cP_{\ga_k}
\stackrel{ \gp_{\ga_k} }{\lra}
\cR_{\ga_k} = \cS_{\ga_k}
\stackrel{\gs_k}{\lra}
\cS_{\ga_0} = \cS_\ga
\stackrel{j}{\lra}
\cM^\cV_{ 1 + \varphi } = \cN$$
where
$\ga_k$ and $\gs_k$ are as in Lemma 3.6
(note that all the models above are weasels).
Since $W$ has the definability property at all ordinals $< \gU$,
the critical points of these two maps are the same.
But the critical point of the latter is $\leq  \pi ( \gk_\ga )$
since we are assuming that $\cS_\ga \not= W$,
and $\crit (i^\cU_{0, \varphi} ) = \crit (E^\cU_\eta )$.
\qed {\tiny \ \ Claim 2}
\vskip 3ex

Pick $t$ as in Claim 1.  Say $t = i^\cU_{0,\eta +1 } (f)(a)$
with $f \in W$, $a \in [ \gn (E^\cU_\eta ) ]^{< \go}$,
and $f: [\crit ( E^\cU_\eta ) ]^{|a|} \lra W$.

\vskip 3ex
\noindent
{\sc Claim 3.}
$\eta = 0$.
\vskip 3ex

\noindent
{\sc Proof.}
Suppose that $\eta > 0$.  Then $\gL_\ga$ is a cardinal in $\cM^\cU_\eta$
and $\lh (E^\cU_\eta ) > \gL_\ga$.
By strong acceptability, for every $\gz < \gL_\ga$,
$$E^\cU_\gg \res ( \gz \cup a ) \in J^{\cM^\cU_\eta}_{\gL_\ga} \ .$$
Also, the $a$-generators of $E^\cU_\eta$ are unbounded in $\gL_\ga$.
Let $G = E^\cU_\eta \res ( \gp (\gk_\ga ) \cup a )$ and let
$A$ be a subset of of $\gp (\gk_\ga )$ which codes $G$.
Then 
$$A \not\in
| \ult (W , G ) |$$
and 
$$ A \in 
| \cM^\cU_{\eta+1} | \ .$$
Let $k$ be the map from $\ult (W , G )$ to $\cM_{\eta +1}^\cU$.
By the choice of $t$,
there is a term $\gt_\psi$ and  parameters 
$b \in [ \gp ( \gk_\ga ) ]^{<\go}$
and 
$c \in [ \mbox{fixed points of }i^\cU_{0,\eta +1}]^{<\go}$,
such that 
for $\gz < \gp (\gk_\ga )$,
$$\gz \in A 
\ \Llra \ 
\gt_\psi^{\cM^\cU_{\eta +1}} [ b , c , t ]
\ \Llra \ 
\gt_\psi^{ult (W , G)}
[b , c , k^{-1}(t)] \ .$$
But then $A \in \ult (W , G )$, a contradiction.
\qed {\tiny \ \ Claim 3}
\vskip 3ex

It remains to see that $\gt (E^\cU_0 )$, the Dodd-projectum
of $E^\cU_0$, is $\leq \gp ( \gk_\ga )$.
It is enough to see that there are no $a$-generators of $E^\cU_0$
$\geq \gp ( \gk_\ga )$.

Suppose, to the contrary,  that 
there is an $a$-generator of $E^\cU_0$
$\geq \gp ( \gk_\ga )$.
Let $G = E^\cU_0 \res ( \gp ( \gk_\ga ) \cup a )$.
Then, from Theorem 2.1.1 it follows that
$G \in J^W_{lh ( E^\cU_0 )}$.
By strong acceptability,
$G \in J^W_{\gL_\ga}$.
But now an argument as in the proof of Claim 3 gives 
a contradiction.  
\qed {\tiny \ \ Lemma 3.11}
\vskip 3ex

The main use of countable closure comes in the proof of the following
lemma.
\vskip 3ex

\noindent
{\sc Lemma 3.13.}
Suppose that $\ga \leq  \gg$ and (3)$_\ga$ holds.
Then (2)$_\ga$ holds.
\vskip 3ex

\noindent
{\sc Proof.}
Suppose (3)$_\ga$ holds and that $\cU$ is a simple iteration tree
of limit length on 
$((W , \cS_\ga ) , \gp ( \gk_\ga ))$.
We shall show that
$\cU$ is well-behaved.

For contradiction, suppose that $\cU$ is ill-behaved.
Take a countable 
elementary submodel of $V_{\gO +1}$
that has as elements $\cU$, $W$, $\cS_\ga$, $\gk_\ga$, and $\gp (\gk_\ga )$.
Let $\psi : M \lra V_{\gO +1}$ be the inverse of the transitive collapse 
of this submodel;
say $\psi (\cU' ) = \cU$, $\psi ( W' ) = W$, $\psi ( \cS' ) = \cS_\ga$,
and $\psi ( \gk' ) = \gp ( \gk_\ga )$.
By the elementarity of $\psi$ and absoluteness as in Section 2 of \cite{St1},
$\cU'$ is a simple, countable, ill-behaved iteration tree on the countable
phalanx $((W' , \cS' ) , \gk' )$.

Recall that $\pi : N \lra V_{\gO +1}$.
The following fact, we probably should have established earlier.

\vskip 3ex
\noindent
{\sc Claim.}
$\cS_\ga = \ult ( \cQ_\ga , E_\gp \res \gp ( \gk_\ga ))$.
\vskip 3ex

\noindent
{\sc Proof.}
The iteration
$$\cR_{\ga_k} = \cS_{\ga_k} \lra \cdots \lra  \cS_{\ga_0} = \cS_\ga$$
is the result of copying the iteration
$$\cP_{\ga_k} = \cQ_{\ga_k} \lra \cdots \lra  \cQ_{\ga_0} = \cQ_\ga$$
using the maps
$\gp_{\ga_k}, \dots , \gp_{\ga_0}= \gp_\ga$.
For $i = k, \dots 0$, let 
$\chi_i : \cQ_{\ga_i} \lra \cS_{\ga_i}$
be the corresponding copying map.
So $\chi_k = \gp_{\ga_k}$ and
by induction using the Shift Lemma,
$\chi_i \res | \cP_{ \ga_i } |
= \gp_{\ga_i}$ for $i = k, \dots 0$.
Set 
$\chi = \chi_0$; then $\chi$ is a map from $\cQ_\ga$ into $\cS_\ga$
and $\chi \res |\cP_\ga | = \gp_\ga$.  Hence
$$E_\chi
= E_\gp \cap ( [ \gp ( \gk_\ga ) ]^{<\go} \times |\cP_\ga | ) \ .$$
So there is a natural $n_\ga$-embedding 
$$h: \ult ( \cQ_\ga , E_\gp \res \gp ( \gk_\ga )) \lra
\cS_\ga$$ such that 
$\crit (h) > \gp ( \gk_\ga )$.
If $\cS_\ga$ is a set premouse, then since $\cS_\ga$ is $\gp (\gk_\ga )$-sound,
equality must hold.
If $\cS_\ga$ is a weasel, then since $\gk (\cS_\ga ) \leq \gp ( \gk_\ga )$
and $c(\cS_\ga ) \in \ran (h)$, and the fixed points of $h$ are thick,
again we may conclude that $h$ is the identity.\\
\qed {\tiny \ \ Claim}
\vskip 3ex

Let $\{ x_0 , \dots , x_m , \dots \}$ and $\{ y_0 , \dots , y_n , \dots \}$
be enumerations of $|W' |$ and $|\cS' |$ respectively.
For each $m<\go$, choose $a_m \in [\gO]^{<\go}$
and $f_m \in N $
such that 
$\psi (x_m ) = \pi ( f_m )(a_m )$.
For each $n<\go$, choose $b_n \in [\gp ( \gk_\ga )]^{<\go}$
and a function $g_n$
such that 
$\dom ( g_n ) = [\gk_\ga ]^{|b_n|}$ and 
$$\psi (y_n ) = [b_n , g_n]^{\cQ_\ga}_{E_\pi \res \gp ( \gk_\ga )} \ .$$
Depending on whether $n_\ga = n(\cQ_\ga , \gk_\ga ) = 0$ or not,
$g_n$ is either an element of $\cQ_\ga$ or is given by a 
$\gS_{n_\ga}$-Skolem term and parameter in $\cQ_\ga$.
If $x_m = y_n < ( ( \gk' )^+ )^{\cS'}$,
then we can and do select $f_m = g_m \in |\cQ_\ga |$ and $a_m = b_n$.
For any $m , n < \go$, put
$c_m = a_0 \cup \cdots \cup a_m$ and 
$d_n = b_0 \cup \cdots \cup b_n$.

Recall the definition, given at the beginning of {\S}2.5,
of the critical point $\gd_a$ 
corresponding to some $a \in [\gO ]^{<\go}$.
Given a first-order formula $\gs ( z_0 , \dots , z_m )$,
let 
$$A_{\gs , m} =
\{ \ u \in [\gd_{c_m}]^{|c_m |} \ |
\ \Wbar \models \gs (f_0^{a_0 , c_m} (u), \dots , 
f_m^{a_m , c_m} (u) ) \ \} \ .$$
Similarly,
for any $\gS_{n_\ga}$-formula $\gt (z_0 , \dots , z_n )$,
let 
$$B_{\gt , n} =
\{ \ u \in [\gd_{d_n}]^{|d_n|} \ |
\ \cQ_\ga \models \gt ( g_0^{b_0 , d_n} (u) , \dots ,
g_n^{b_n , d_n} (u) ) \ \} \ .$$

Since $ ^{\go} (X \cap V_\gO ) \subseteq X 
\buildrel \gp
\over \prec 
V_{\gO + 1}$, there are order preserving maps
$s : \bigcup_{m<\go} a_m \lra \gObar$ and 
$t : \bigcup_{n<\go} b_n \lra \gk_\ga$
such that
$$A_{\gs , m} \in (E_\gp )_{c_m} \ \Lra \ s \, '' \, c_m \in A_{\gs , m}\ ,$$
$$B_{\gt , n} \in (E_\gp )_{d_n} \ \Lra \ t \, '' \, d_n \in B_{\gt , n}\ ,$$
and
$$x_m = y_n < ((\gk' )^+ ))^{\cS'} 
\ \Lra \ s \, '' \, a_m = t\, '' \, b_n \ .$$
Define
$\varphi_0 : W' \lra \Wbar$ and $\varphi_1 : \cS' \lra \cQ_\ga$
by
$$\varphi_0 (x_m ) = f_m ( s\, '' \, a_m )$$
and 
$$\varphi_1 (y_n ) = g_n ( t\, '' \, b_n )$$
Then
$$\varphi_0 \res ((\gk' )^+ )^{\cS'} =
\varphi_1 \res ((\gk' )^+ )^{\cS'} \ .$$
Using the {\L}o{\'s} Theorem,
we see that $\varphi_0$ is an elementary embedding of $W'$ into $\Wbar$,
and $\varphi_1$ is a weak $n_\ga$-embedding 
of $\cS'$ into $\cQ_\ga$.  Let $\chi : \cQ_\ga \lra \cS_\ga$ be 
the natural $n_{\ga}$-embedding, as in the proof of the Claim.  
We also know that 
$\chi \circ \varphi_1 = \psi \res | \cS'|$.
Since $\chi$ is $\gS_{n_\ga + 1}$-elementary,
$\varphi_1$ is $\gS_{n_\ga}$-elementary,
and $\psi \res | \cS'|$ is fully elementary,
it follows that, in fact,
$\varphi_1$ is $\gS_{n_\ga + 1}$-elementary.
We also know that 
$\pi \circ \varphi_0 = \psi \res |W' |$.

Now copy $\cU'$ 
using the pair of maps $(\varphi_0 , \varphi_1 )$;
the result is
an ill-behaved iteration tree on $((\Wbar , \cQ_\ga ), \gk_\ga )$.
This contradicts (3)$_\ga$.
\qed {\tiny \ \ Lemma 3.13}
\vskip 3ex

We next work towards proving (3)$_\ga$ and (1)$_\ga$
using instances of (4)$_\gb$.

\vskip 3ex
\noindent
{\sc Lemma 3.14.}
Suppose that $\ga \leq \gg$, (1)$_\gb$ holds for every $\gb < \ga$,
and
(4)$_\ga$ holds.  Then
there is an iteration tree $\cTtil$ on $W$
and an elementary embedding $k : \Wbar \lra \cN$
such that:
\begin{list}{}{}
\item[(a)]{$\cTtil$ extends $\cT \res ( \eta ( \ga ) + 1 )$.}
\item[(b)]{ $\cTtil$ has a final model and $\cN$ is an initial segment 
of the final model of $\cTtil$.}
\item[(c)]{$\lh (E^\cTtil_\eta ) \geq \gl_\ga$ whenever 
$\eta ( \ga ) \leq \eta < \lh ( \cTtil )$.}
\item[(d)]{
$\crit (k) \geq \gk_\ga$.}
\end{list}
\vskip 3ex

\noindent
{\sc Proof.}
Since (1)$_\gb$ holds for every $\gb < \ga$,
$\cP_\ga$ and $\Wbar$ agree below $\gl_\ga$.
We must allow for the possibility, though, that
they disagree at $\gl_\ga$.
Both $(\cPvec \res (\ga + 1 ) , \glvec \res \ga )$ and 
and $( ( \cPvec \res \ga , \Wbar ) , \glvec \res \ga )$
are iterable, very  special phalanxes of premice, 
since we are assuming (4)$_\ga$.
Compare these two phalanxes.  This results in iteration trees
$\cU$ on 
$(\cPvec \res (\ga + 1 ) , \glvec \res \ga )$
and $\cV$ on
$( ( \cPvec \res \ga , \Wbar ) , \glvec \res \ga )$
with $\lh ( \cU ) = \ga + \varphi +1 = \lh ( \cV )$
for some $\varphi$,
and either $\cM^\cU_{\ga + \varphi}$
is and initial segment of $\cM^\cV_{\ga + \varphi}$
or vice-versa.
All extenders used on $\cU$ and $\cV$
have length $\geq \gl_\ga$.
The iteration tree $\cTtil$ whose existence 
is asserted by Lemma 3.14 is just $\cU$ reorganized in the obvious
way.  We shall show that $\ga = \rt^V (\ga + \varphi ) $ and that 
$[\ga , \ga + \varphi ]_U \cap \cD^\cU = \emptyset$.  We will
also show that $\cM^\cV_{\ga + \varphi}$ is an initial segment
of $\cM^\cU_{\ga + \varphi}$.  Once this is accomplished,
we can set $\cN = \cM^\cV_{\ga + \varphi}$ and 
$k  = i^\cV_{\ga , \ga + \varphi}$, and Lemma 3.14 follows.

\vskip 3ex
\noindent
{\sc Claim 1.}
Let 
$$\etatil = 
\sup 
\left( \{ \eta < \eta ( \ga ) \ | \ E_\eta \not= \emptyset \} \right) \ .$$
Then $\Wbar = \Wbar_\etatil$ and $W_\etatil = W_{\eta ( \ga )}$.
If $\Ebar_\etatil \not= \emptyset$ (that is, if we are not padding
at stage $\etatil$ in the construction of $\cTbar$), then 
$E^\cV_\ga = \Ebar_\etatil$.
If $\etatil = \eta ( \ga )$, then
$E^\cU_\ga = E_{\eta ( \ga )}$.
\vskip 3ex

\noindent
{\sc Proof.}
If $\etatil < \eta ( \ga )$, then $\etatil T \eta ( \ga )$
and $(\etatil , \eta ( \ga ) ]_T$ consists only of padding.
That is, $W_\etatil = W_{\eta ( \ga )}$.

Suppose that $\eta < \etatil$ and $\Ebar_\etatil \not= \emptyset$.
as we noted,
$\lh ( \Ebar_\eta ) \geq \gl_\ga$.
This implies that $\gn ( E_\xi ) \geq \gl_\ga $
whenever $\xi > \eta $ and $E_\xi \not= \emptyset$.
Thus $E_\xi = \emptyset$ whenever $\eta < \xi < \eta ( \ga )$,
so $\etatil \leq \eta$, a contradiction.  Therefore 
$\Wbar_\etatil = \Wbar$.

If $\Ebar_\etatil \not= \emptyset$, then $\Ebar_\etatil$
is part of the least disagreement between $\Wbar_\etatil$ and
$W_\etatil$, hence between $\Wbar$ and $W_{\eta ( \ga )}$,
hence between $\Wbar$ and $\cP_\ga$.

If $\etatil = \eta ( \ga )$,
then $E_{\eta ( \ga ) }$ is part of the least disagreement
between $\Wbar_\etatil$ and $W_{\eta ( \ga )}$,
hence between $\Wbar$ and $\cP_\ga$.
\qed {\tiny \ \ Claim 1}
\vskip 3ex

Say $\gb_0 = \rt^U (\ga + \varphi )$
and $\gb_1 = \rt^V (\ga + \varphi )$.
Then $\gb_0$, $\gb_1 \leq \ga$.

\vskip 3ex
\noindent
{\sc Claim 2.} 
$\gb_1  = \ga$, that is, the last model on $\cV$ follows 
$\Wbar$.
\vskip 3ex

Assume to the contrary that $\gb_1 < \ga$.

\vskip 3ex
\noindent
{\sc Subclaim A.}
$\cM^\cU_{\ga + \varphi}$ is an initial segment of $\cM^\cV_{\ga + \varphi}$.
\vskip 3ex

\noindent
{\sc Proof.}
Either because $\cP_{\gb_1}$ is $\gk_{\gb_1}$-sound and 
$\cM^\cV_{\ga + \varphi }$ is not,
or because $\cP_{\gb_1}$ is a universal weasel.
\qed {\tiny \ \ Subclaim A}
\vskip 3ex

\noindent
{\sc Subclaim B.}
$\cM^\cU_{\ga + \varphi}
= \cM^\cV_{\ga + \varphi}$.
\vskip 3ex

\noindent
{\sc Proof.}
This is clear unless $\gb_0 = \ga$ and $\cM^\cU_{\ga + \varphi} = \cP_\ga$.
So suppose that $\cP_\ga $ is a proper initial segment of 
$\cM^\cV_{\ga + \varphi}$.
Then $\cP_\ga$ is a sound, set premouse that projects to $\gk_\ga$,
so $$\OR^{\cP_\ga} < ( \gk_\ga^+ )^{\cM^\cV_{\ga + \varphi}}$$
But by Claim 1, the right hand side is $\gl_\ga$.
Thus
$$\OR^{\cP_\ga} \leq
( \gk_\ga^+ )^{\Wbar_\theta} \leq (\gk_\ga^+ )^{W_\theta}
<
(\gk_\ga^+ )^{W_{\eta ( \ga ) + 1}}
\leq
\lh (E_{\eta ( \ga ) }) \leq \OR^{\cP_\ga} \ ,$$
which is absurd.
\qed {\tiny \ \ Subclaim B}
\vskip 3ex

Put $\cM = \cM^\cU_{\ga + \varphi}
= \cM^\cV_{\ga + \varphi}$.
Let  $\gs = \max ( \ [ 0,\eta ( \gb_0 ) ]_T \ \cap \ 
[ 0,\eta ( \gb_1 ) ]_T \ )$.  We have two different ways of iterating
$W$ to $\cM$:
$$\begin{array}{ccccccc}
\ & \ & \ & \ & \! W_{\eta ( \gb_0 )} \ \supseteq \cP_{\gb_0} & \lra &
\cM^\cU_{\ga + \varphi} = \cM \cr
\ & \  & \ & \! \nearrow & \ & \ & \  \cr
W & \lra & W_\gs & \ & \ & \ & \  \cr
\ & \ & \ & \! \searrow & \ & \ & \  \cr
\ & \ & \ & \ & \! W_{\eta ( \gb_1 )} \ \supseteq \cP_{\gb_1} & \lra &
\cM^\cV_{\ga + \varphi} = \cM \cr
\end{array}$$
The arrows do not necessarily
indicate maps, as we may have {\it drops} across the top or bottom.
For the purposes of this proof only,
by a drop across the top we mean one of the following:
\begin{list}{}{}
\item[I.]{an ordinal $\eta + 1 \leq_T \eta (\gb_0 )$
such that $\eta + 1$ is a drop of any kind in the sense of $\cT$, or}
\item[II.]{the ordinal $\ga + \xi + 1$ 
such that 
$\ga + \xi +1$ is least in 
$(\gb_0 , \ga + \varphi ]_U$, on the condition that
$\cP_{\gb_0}$ is a proper initial  
segment of $W_{\eta ( \gb_0 )}$, or}
\item[III.]{an ordinal $\ga + \xi + 1$ that is a drop of any kind in
the sense of $\cU$, so long as
$\ga + \xi + 1 \leq_U \ga + \varphi$ and 
$\gb_0 \not= \pred^\cU ( \ga + \xi + 1 )$.}
\end{list}
And replacing
$\gb_0$ and $\cU$ by $\gb_1$ and $\cV$ defines a drop across the bottom.
There are many cases, depending on whether or not there
are drops, and, if there are, where the 
last drops across the top and bottom
occur.  In all these cases, the contradiction
is roughly like in the Comparison Lemma 7.1 of \cite{MiSt}.
We shall give the argument in a single illustrative case,
and leave the other cases to the reader.

\vskip 3ex
\noindent
{\sc Illustrative Case.}
Suppose that the last drop across the top is $\eta + 1$,
a drop by condition I, and the
last drop across the bottom is $\ga + \xi + 1$, a drop by
condition II or III.
Let $G_0 = E_\eta$ and $G_1 = E^\cV_{\ga + \xi}$.
If $\gn = \min ( \gn ( G_0 ), \gn ( G_1 ))$,
then by a standard argument, $G_0 \res \gn = G_1 \res \gn$.
Since $\eta < \eta ( \gb_0 ) \leq \eta ( \ga )$,
$\gn ( G_0 ) \leq \gk_\ga$.
$\lh (G_0 ) \leq \gl_\ga$, 
as otherwise, $\lh (G_0 ) > \gl_\ga$,
so $\gl_\ga$ is a cardinal in $J^{W_\eta}_{lh (G_0 )}$,
so $\gn (E_\eta ) \geq \gl_\ga$, so $\eta \geq \eta ( \ga )$,
but it is not.
If $\gn ( G_0 ) = \gk_\ga$,
then $\lh ( G_0 ) = \gl_\ga$, as there are no cardinals
between $\gk_\ga$ and $\gl_\ga$ in $W_\theta$.
Therefore $$\lh ( G_0 ) \leq \gl_\ga \leq \lh (G_1 )$$
and $$\gn ( G_0 ) \leq \gk_\ga \leq \gn (G_1 ) \ .$$
Since $\eta < \eta ( \ga )$,
$G_0 \not\in | W_{\eta ( \ga )} |$, so $G_0 \not\in |\cP_\ga |$.
Suppose first that $\gn ( G_0 ) < \gn ( G_1 )$.
Then, by the initial segment condition on $\cM^\cV_{\ga + \xi \ }$,
$G_0$ is an element of $\cM^\cV_{\ga + \xi}$ constructed at a level before
$$\gl_\ga = (\gk_\ga^+ )^{ J^{\cM^\cV_{\ga + \xi}}_{lh (G_1 )} } \ .$$
But then $$G_0 \in J^\Wbar_{\gl_\ga}
=J^{\cP_\ga}_{\gl_\ga} \ ,$$
a contradiction.
So $\gn (G_0 ) = \gn (G_1 ) = \gk_\ga$.
But then $G_0 = G_1$, and this is an extender of length $\gl_\ga$.
So $\xi = 0$ and $G_1$ is an extender on the $\Wbar$-sequence.
Note that $\eta < \etatil$, so $\Wbar_\eta = \Wbar$ by Claim 1.
Therefore $$E_\eta = E^{W_\eta}_{\gl_\ga} = G_0 = G_1 =
E^{\Wbar_\eta}_{\gl_\ga} = \Ebar_\eta$$
which contradicts the rules for comparison by which we formed $\cT$ and
$\cTbar$.\\
\qed {\tiny \ \ Illustrative Case} \qed {\tiny \ \ Claim 2}
\vskip 3ex

By an almost identical argument, we see that there is no dropping
along the main branch of $\cV$.  Again, we leave the details to the 
reader.

\vskip 3ex
\noindent
{\sc Claim 3.} $[\ga , \ga + \varphi ]_U \cap \cD^\cU = \emptyset$.
\vskip 3ex

\noindent
{\sc Claim 4.} $\cM^\cV_{\ga + \varphi}$ is an initial segment of
$\cM^\cU_{\ga + \varphi}$.
\vskip 3ex

\noindent
{\sc Proof.}
This is clear unless $\cM^\cU_{\ga + \varphi} = \cP_\ga$.
But we can show that $\cP_\ga$ is not a proper initial segment
of $\cM^\cV_{\ga + \varphi}$ just as in the proof of Subclaim A.
\qed {\tiny \ \ Claim 4} \qed {\tiny \ \ Lemma 3.14}
\vskip 3ex

\noindent
{\sc Lemma 3.15.} 
Suppose that $\ga \leq \gg$, (1)$_\gb$ holds for every $\gb < \ga$,
and (4)$_\ga$ holds.  Then (1)$_\ga$ holds.
\vskip 3ex

\noindent
{\sc Proof.}
Suppose that Lemma 3.15 fails.  Then there is an $\eta$ such that
$\Ebar_\eta \not= \emptyset$ and $\lh (\Ebar_\eta ) \leq \gl_\ga$.
Since (1)$_\gb$ holds for every $\gb < \ga$,
$\lh (\Ebar_\eta ) = \gl_\ga$.  Therefore 
$$\crit (\Ebar_\eta ) 
< \gk_\ga \leq \gn (\Ebar_\eta ) < \gl_\ga = \lh (\Ebar_\eta ).$$

Now continue the comparison begun in the proof of Lemma 3.14.
Then $\eta \leq \eta ( \ga )$,
$E^\cV_\ga = \Ebar_\eta$,
and $\pred^\cV ( \ga +1 ) < \ga = \rt^V (\ga + \varphi )$.
Let $\ga + \xi +1$ be least in $(\ga , \ga + \varphi ]_V$.
Then $\gl_\gb \leq \crit (E^\cV_{\ga + \xi }) < \gn (\Ebar_\eta )$
for every $\gb < \ga$, which can mean only that 
$\crit (E^\cV_{\ga + \xi}) = \gk_\ga$.

Since we hit $\Ebar_\eta$, an extender of length $\gl_\ga$, 
at the first step in building $\cV$,
$$\gl_\ga = (\gk_\ga^+ )^{\cM^\cV_{\ga + \xi }} < (\gk_\ga^+ )^\Wbar.$$
So $\ga + \xi + 1 \in \cD^\cU$, contradicting Lemma 3.14.
\qed {\tiny \ \ Lemma 3.15}
\vskip 3ex

\noindent
{\sc Lemma 3.16.}
Suppose that (4)$_\gb$ holds for every $\gb < \ga$.
Then (3)$_\ga$ holds.
\vskip 3ex

\noindent
{\sc Proof.}
By Lemma 3.15, we know that (1)$_\gb$ holds for every $\gb < \ga$.
So by Lemma 3.14, for every $\gb < \ga$, 
there is an iteration tree $\cTtil_\gb$ on $W$ with $\lh ( \cTtil _\gb )
= \varphi_\gb + 1$ such that $\cTtil_\gb$ extends
$\cT \res (\eta ( \gb ) +1 )$ and there is an elementary embedding $k_\gb$
with critical point at least $\gk_\gb$
from $\Wbar$ into an initial segment $\cN_\gb$ of the last model
$\cM^{\cTtil_\gb}_{\varphi_\gb}$ of $\cTtil_\gb$.
All extenders $E^{\cTtil_\gb}_\eta$ have length at least $\gl_\gb$
whenever $\eta ( \gb ) \leq \eta < \lh ( \cTtil_\gb )$.
Let $k_\ga = \id \res |\cQ_\ga |$,
$\cNvec = \langle \ \cN_\gb \ | \ \gb < \ga \ \rangle$, and
$\kvec = \langle \ k_\gb \ | \ \gb \leq \ga \ \rangle$.
The proof of the following fact is 
is similar to the proof of the main result in Section 9 of \cite{St1};
we omit the details.

\vskip 3ex
\noindent
{\sc Fact 3.16.1.}
$( ( \cNvec , \cQ_\ga ) , \gkvec \res \ga )$ is an iterable phalanx.
\vskip 3ex

Suppose that $\cU$ is an iteration tree on 
$( (\Wbar , \cQ_\ga ) , \gk_\ga )$.
We wish to copy $\cU$ to an iteration tree $\cV$ on 
$( ( \cNvec , \cQ_\ga ) , \gkvec \res \ga )$
using the system $\kvec$.  
By Fact 3.16.1, this is enough.
Note that there are two base models 
on $\cU$, while there will be
$\ga$-many base models on $\cV$.

To avoid some triple superscripts,
we shall use the notation:
$$\cM^* (\cU , 1 + \eta + 1) = (\cM^*_{1 + \eta + 1} )^\cU$$ 
and
$$\cM^* (\cV , \ga + \eta + 1) = (\cM^*_{\ga + \eta + 1} )^\cV$$ 
in the discussion immediately below.

There is really only one subtlety in the copying construction giving $\cV$.
Suppose that for some $1 + \eta  + 1 < \lh ( \cU )$,
$E^\cU_{1 + \eta}$ has critical point $\gk_\gb$ for some $\gb < \ga$.
Then $E^\cV_{\ga + \eta} = k_{1 + \eta} ( E^\cU_{1 + \eta} )$,
with the usual understanding that  
$k_{1 + \eta} ( \Fdot^{ \cM^\cU_{1 + \eta} } ) =
\Fdot^{ \cM^\cV_{\ga + \eta} }$.
Since 
$$k_{1 + \eta} \res \gk_\ga = k_\ga \res \gk_\ga = \id \res \gk_\ga \ ,$$
$\crit ( E^\cV_{\ga + \eta} ) = \gk_\gb$.
But it is possible that $\crit (k_\gb ) = \gk_\gb$, in which case
the Shift Lemma 5.2 of \cite{MiSt} does not literally apply.
But the proof works if we copy as follows.
In this situation, $( \cM^*_{1 + \eta + 1} )^\cU$ is an initial segment
of $\Wbar$
and $( \cM^*_{\ga  + \eta + 1} )^\cV$ is an initial segment
of $\cN_\gb$.
Given coordinates $a \in [ \lh (E^\cU_{1 + \eta} ) ]^{< \go}$
and a function 
$$f : 
[ \crit (E^\cU_{1 + \eta} ) ]^{|a|} \lra | ( \cM^*_{1 + \eta + 1} )^\cU |$$
with either $f \in |  ( \cM^*_{1 + \eta + 1} )^\cU |$ or given by 
a Skolem term and parameter in $ ( \cM^*_{1 + \eta + 1} )^\cU$,
let 
$$k_{1 + \eta + 1} \left( \ 
\left[ a , f \right]^{\cM^* (\cU , 1 + \eta + 1)}_{E^\cU_{1 + \eta}} \ 
\right)
= \left[ \ k_{1+ \eta} (a) \ , 
\ k_\gb ( f ) \res [ \gk_\gb ]^{|a|} 
\ \right]^{
\cM^* (\cV , \ga + \eta + 1)}_{ E^\cV_{\ga + \eta} } \ .$$
The proof of the Shift Lemma still works.  For example,
given 
$$x 
= \left[ a,f \right]^{\cM^* (\cU , 1 + \eta + 1 )}_{ E^\cU_{1 + \eta} }$$
and 
$$y 
= \left[ a,g \right]^{\cM^* (\cU , 1 + \eta + 1)}_{ E^\cU_{1 + \eta} },$$
if
$X = \{ \ u \in [\gk_\gb ]^{|a|} \ | \ f(u)=g(u) \ \}$,
then
$x = y$
if and only if
$X \in ( E^\cU_{1 + \eta} )_a$
if and only if
$$k_{1 + \eta} (X) = X = k_\gb ( X ) 
\cap \gk_\gb \in ( E^\cV_{1 + \eta} )_{k_{1 + \eta}(a)}$$
if and only if
$k_{1 + \eta + 1} (x) = 
k_{1 + \eta + 1} (y)$.

The rest of the details are left to the reader. Note that
the above variation on the usual copying construction is used
in Section 6 of \cite{St1}.
\qed {\tiny \ \ Lemma 3.16}
\vskip 3ex

\noindent
{\sc Lemma 3.17.}
Suppose that (5)$_\ga$ holds.  Then (4)$_\ga$ holds.
\vskip 3ex

\noindent
{\sc Proof.}
We can use the system of maps 
$\gsvec = 
(\langle \ \gp_\gb \ | \ \gb < \ga \ \rangle), \gp)$
to copy a given iteration tree $\cU$ on the phalanx
$( (\cPvec \res \ga , \Wbar ) , \glvec \res \ga )$ 
to an iteration tree $\gsvec '' \cU$
on the phalanx
$( (\cRvec \res \ga , W ) , \gLvec \res \ga )$.
Notice that $\gsvec '' \cU$ is a special iteration tree,
since every extender used on it has critical point in the range
of $\gp$.  So
$\gsvec '' \cU$ is well-behaved.  The usual argument now
shows that $\cU$ is well-behaved.
\qed {\tiny \ \ Lemma 3.17}
\vskip 3ex

\noindent
{\sc Lemma 3.18.}
Suppose that (6)$_\ga$ holds.  Then (5)$_\ga$ holds.
\vskip 3ex

\noindent
{\sc Proof.}
Let $\cU$ be a special iteration tree on
$((\cRvec , W), \gLvec \res \ga )$.
We wish to {\it enlarge} $\cU$ to a special iteration tree $\cV$
on $((\cSvec , W), \gLvec \res \ga )$.
This will differ significantly from constructions that
we are familiar with.  For example, although we will have
$\lh ( \cU ) = \lh ( \cV )$, it will be possible that
the tree structures
$U$ and $V$ are different.
Also,
it will be possible 
that for some $\ga + \eta + 1 < \lh ( \cU )$,
$E^\cU_{\ga + \eta} \not= \emptyset$ while 
$E^\cV_{\ga + \eta} = \emptyset$.
Before we describe the construction of $\cV$,
let us make some observations about $\cU$.
Suppose that $\ga \leq \ga + \eta + 1 < \lh ( \cU )$.
Let $\gb = \rt^\cU ( \ga + \eta )$.  We identify two cases:

\vskip 3ex
\noindent
{\sc Case A.}
$\Fdot^{\cM^\cU_{\ga + \eta}}$
is not an extender over $\cM^\cU_{\ga + \eta}$.
\vskip 3ex

\noindent
{\sc Case B.}
Otherwise.
\vskip 3ex

Case A occurs exactly when the following three conditions hold:
\begin{list}{}{}
\item[1.]{$\gb < \ga$}
\item[2.]{$\cR_\gb$ is not a premouse; equivalently, $\cR_\gb \not=
\cS_\gb$}
\item[3.]{$(\gb , \ga + \eta ]_U \cap \cD^\cU = \emptyset$}
\end{list}
In Case A, we also have that:
\begin{list}{}{}
\item[4.]{$\cR_\gb$ is an active protomouse that is active of type
	I or II}
\item[5.]{$m_\gb = 0$}
\item[6.]{$\deg^\cU ( \ga + \eta ) = 0$}
\item[7.]{There is a unique $\gbtil < \gb$ such that
	$\gmdot^{\cR_\gb} = \gp (\gk_{\gbtil} )$}
\item[8.]{$\crit ( i^\cU_{\gb , \ga + \eta} ) 
	> ( \gmdot^+ )^{\cR_\gb}$}
\item[9.]{$\gmdot^{\cM^\cU_{\ga + \eta}}
	= \gmdot^{\cR_\gb} $}
\item[10.]{$\Fdot^{\cR_\gb}$ is an extender over $\cR_{\gbtil}$}
\end{list}
Moreover, the relationship between $\cR_\gb$ and $\cS_\gb$ is as described
in Lemmas 3.6 and 3.7.

Note that in Case A, if $E^\cU_{\ga + \eta} = \Fdot^{\cM^\cU_{\ga + \eta}}$,
then $\pred^\cU  ( \ga + \eta + 1) = \gbtil$,
$$\ga + \eta + 1 \not\in \cD^\cU \ ,$$
and $\deg^\cU ( \ga + \eta + 1 ) = m_{\gbtil}$.
In this case,
$E^\cU_{\ga + \eta}$
measures precisely those
subsets of 
its critical point that are 
in $\cR_\gbtil = ( \cM^*_{\ga + \eta + 1} )^\cU$.
Note that in other cases,
it is possible that $E^\cU_{\ga + \eta}$ measures more subsets 
of its critical point than are present in $( \cM^*_{\ga + \eta + 1} )^\cU$.
This happens, for example, if $E^\cU_{\ga + \eta}$ is an extender
over $$J^{\cM^\cU_{\ga + \eta}}_{lh ( E^\cU_{\ga + \eta} )}$$ and
$\crit ( E^\cU_{\ga + \eta} ) = \gp ( \gk_\gb )$ for some $\gb < \ga$
such that $\cR_\gb$ is not a premouse.  In such cases, we are using
the full force of Definition 2.4.2.

Consider any $\gb < \ga$.
If $\cR_\gb = \cS_\gb$, then let $\cL_\gb$ be this premouse.
If, on the other hand, $\cR_\gb \not= \cS_\gb$,
then let $\cL_\gb = \decap ( \cR_\gb )$.
Let $\cL_\ga = W$.
For all $\gb \leq \ga$, let $e_\ga$
be the identity map on $| \cL_\gb |$.
These are the starting maps for our 
enlargement.  
Either $\cL_\gb$ equals $\cS_\gb$, or else
$\cL_\gb$ is $\cS_\gb$ cut off at a successor cardinal
of $\cS_\gb$.

Suppose, for the moment, that $\cR_\gb \not= \cS_\gb$.
Let 
$$\cR_{\gb_k} = \cS_{\gb_k} \lra \cdots \lra \cS_{\gb_0} = \cS_\gb$$
be the decomposition of $\cS_\gb$ as in Lemma 3.6,
with $\gs_i : \cS_{\gb_i} \lra \cS_\gb$ for $i = 0,1, ... , k$
the natural maps.  Note that $n_\gb = n_{\gb_0} = n_{\gb_1}
= \cdots = n_{\gb_k} \leq \go$; call this ordinal $n$.
Also suppose that $n < \go$.
Then
$$ \gr_{n +1}^{\cS_\gb} \leq \gp ( \gk_\gb ) < \gr_n ( \cS_{\gb_i} )$$
and for $i = 0, 1, \dots , k$,
$$ \gr_n ( \cS_\gb ) = \gs_i ( \gr_n ( \cS_{\gb_i} ))
> \gs_i ( \gp ( \gk_{\gb_i} )) \ .$$
Also,
for $i = 1, \dots , k$,
$$\gs_i (\gp (\gk_{\gb_i} )) > 
\gs_{i-1} ( \OR \cap \cR_{ \gb_{i-1} } )
= \OR \cap \gs_{i-1} ( \cL_{ \gb_{i-1} } ) \ .$$
In particular,
$$\gr_{n+1}^{\cS_\gb} < \gp ( \gk_\gb ) < \OR^{\cL_\gb} 
< \gr_n ( \cS_\gb )$$
whenever $\cL_\gb$ is a proper initial segment of $\cS_\gb$.

By induction on $\ga + \eta < \lh ( \cU )$, we wish to define the initial
segments $\cV \res ( \ga + \eta ) $ of a special iteration tree
$\cV$ on 
$((\cSvec \res \ga , W ) , \gLvec \res \ga )$ together with 
initial segments $\cL_{\ga + \eta}$ of $\cM^\cV_{\ga + \eta}$
and maps $e_{\ga + \eta} : \cM^\cU_{\ga + \eta} \lra \cL_{\ga + \eta}$.
The maps will satisfy certain inductively maintained agreement, elementarity,
and commutativity properties.  In order to state these properties,
we first make the following definition.

\vskip 3ex
\noindent
{\sc Remark:}
There is a slight lie here.
$\cV$ will involve a generalization
of padding, to be explained at the beginning of Case 4 below.
\vskip 3ex

Suppose that 
$\ga + \eta < \lh (\cU )$ and that $\gb = \rt^\cU (\ga + \eta )$.
There is a unique pair of (finite) sequences $\gb_0 > \cdots > \gb_\ell$
and $\eta_0 < \cdots < \eta_{\ell - 1} < \eta_\ell$ that are of maximal
length with the following properties.
\begin{itemize}
\item	$\eta_\ell = \eta$.
\item	If $0 \leq i \leq \ell$,
	then $\gb_i = \rt^\cU ( \ga + \eta_i )$; so $\gb_\ell = \gb$.
\item	If $0 < i \leq \ell$,
	then $\gb_i = \pred^\cU (\ga + \eta_{i-1} + 1 )
	\leq_U (\ga + \eta_i )$.
\item	If $0 < i \leq \ell$,
	then $\ga + \eta_{i-1}$ falls under Case A
	and $E^\cU_{\ga + \eta_{i-1}} 
	= \Fdot^{\cM^\cU_{\ga + \eta_{i-1}}}$.
\end{itemize}
We shall call $( \gbvec , \etavec)$ the {\bf trace} of $\ga + \eta$.
It is not difficult to see that the decomposition of $\cS_{\gb_0}$
in the sense of Lemma 3.6 is
$$\cR_{\gb_k} 
= \cS_{\gb_k} \lra \cdots \lra \cS_{\gb_\ell} \lra \cdots \lra \cS_{\gb_0}$$
for some $k \geq \ell$.
The following diagram may serve as a useful reference in what is to
come.

\begin{picture}(400,330)(0,0)
\put(0,25){\underline{{\sc Figure 1}}}
\put(100,25){$\cR_{\gb_2}$}
\put(150,50){$\cR_{\gb_1}$}
\put(200,75){$\cR_{\gb_0}$}

\put(250,25){$\cS_{\gb_2}$}
\put(300,50){$\cS_{\gb_1}$}
\put(350,75){$\cS_{\gb_0}$}

\put(270,38){\vector(2,1){17}}
\put(320,60){\vector(2,1){17}}
\put(333,56){$\gs_1$}
\put(281,29){\tiny{$\Fdot^{\cR_{\gb_1}}$}}

\put(135,27){\vector(1,0){92}}
\put(185,52){\vector(1,0){92}}
\put(235,77){\vector(1,0){92}}

\put(202,100){\vector(0,1){10}}
\put(202,130){\vector(0,1){10}}

\put(352,100){\vector(0,1){10}}
\put(352,130){\vector(0,1){10}}

\put(195,150){$\cM^\cU_{\ga + {\eta_0}}$}
\put(345,150){$\cM^\cV_{\ga + {\eta_0}}$}
\put(235,152){\vector(1,0){92}}
\put(260,160){$e_{\ga + \eta_0}$}

\put(350,173){$\|$}

\put(345,193){$\cM^\cV_{\ga + {\eta_0} +1}$}
\put(145,193){$\cM^\cU_{\ga + {\eta_0} +1}$}

\put(195,198){\vector(1,0){132}}
\put(240,205){$e_{\ga + \eta_0 + 1}$}
\put(155,75){\vector(0,1){100}}
\put(114,130){$\Fdot^{\cM^\cU_{\ga + \eta_0}}$}

\put(155,210){\vector(0,1){10}}
\put(155,230){\vector(0,1){10}}
\put(350,210){\vector(0,1){10}}

\put(350,210){\vector(0,1){10}}
\put(350,230){\vector(0,1){10}}

\put(345,250){$\cM^\cV_{\ga + \eta_1}$}
\put(145,250){$\cM^\cU_{\ga + \eta_1}$}

\put(195,250){\vector(1,0){132}}
\put(240,260){$e_{\ga + \eta_1}$}

\put(90,295){$\cM^\cU_{\ga + \eta_1 + 1}$}
\put(345,295){$\cM^\cV_{\ga + \eta_1 + 1}$}

\put(350,273){$\|$}

\put(150,297){\vector(1,0){177}}
\put(235,310){$e_{\ga + \eta_1 + 1}$}

\put(103,60){\vector(0,1){215}}
\put(50,200){$\Fdot^{\cM^\cU_{\ga + \eta_1}}$ }

\end{picture}

Let $\ga + \eta < \lh ( \cU )$ be given, and suppose it has trace
$(\gbvec , \etavec )$.  Put $\gb = \rt^\cU ( \ga + \eta ) = \gb_\ell$.
We inductively maintain the following properties.

\begin{list}{}{}
\item[
\underline{The tree, drop, and degree structure of $\cV$}]{\ \\

For any $\xi \leq \ga + \eta$, 
$\xi \leq_V (\ga + \eta)$ if and only if
$\xi \leq_ U (\ga + \eta_i )$ for some $i \in \{ 0 , \dots , \ell \}$.
This determines the tree structure $V$ of $\cV$.\\

$\cD^\cV \cap (\ga + \eta + 1) = \cD^\cU \cap (\ga + \eta + 1)$,
so in the end, we shall have that the drop structures for $\cU$ and
$\cV$ are the same.\\

Suppose that $\ga + \eta$ falls under Case A. Then
$$\deg^\cV ( \ga + \eta ) = n_{\gb_0} = \cdots = n_{\gb_\ell} =
n_\gb..$$
Recall that $\deg^\cU ( \ga + \eta ) = 0$ in Case A;
so it is possible that 
$\deg^\cV (\ga + \eta)$ is different from
$\deg^\cV (\ga + \eta)$.\\

On the other hand, suppose that $\ga + \eta$ falls under Case B.  Then
$\deg^\cV ( \ga + \eta ) = 
\deg^\cU ( \ga + \eta )$.\\}

\item[
\underline{The premouse $\cL_{\ga + \eta}$}]{\ \\

Suppose that $\ga + \eta$ falls under Case A.
Then $\cL_{\ga + \eta}$ 
is a proper initial segment of $\cM^\cV_{\ga + \eta}$
and $\OR^{\cL_{\ga + \eta}}$ is a successor cardinal in $\cM^\cV_{\ga + \eta}$.
If $\cM^\cV_{\ga + \eta}$ is a set premouse, then
$$\gr^{\cM^\cV_{\ga + \eta}}_{n_{\gb_0} +1} \leq
\gp ( \gk_{\gb_0} ) < \OR^{\cL_{\ga + \eta}}
< \gr_{n_{\gb_0}} ( \cM^\cV_{\ga + \eta} ) \ .$$
More precisely,
$$\cL_{\ga + \eta} = ( i^\cV_{\gb_0 , \ga + \eta} \circ \gs_\ell )
( \cL_\gb )$$
where, as before, $\gs_\ell : \cS_{\gb_\ell} \lra \cS_{\gb_0}$
is the $\ell$'th decomposition map for $\cS_{\gb_0}$.\\

On the other hand, suppose that $\ga + \eta$ falls under Case B.
Then $\cL_{\ga + \eta} = \cM^\cV_{\ga + \eta}$.\\}

\item[
\underline{Agreement}]{\ \\

Suppose that $\ga + \xi \leq_U \ga + \eta_i$ for some 
$i \in \{ 0 , \dots , \ell \}$.
Then $$e_{\ga + \eta} \res \gn^\cU_\xi
= e_{\ga + \xi} \res \gn^\cU_\xi.$$\\}

\item[\underline{Elementarity}]{\ \\

Suppose that $\ga + \eta$ falls under Case A.
Then
$e_{\ga + \eta}$ is an elementary embedding of 
$\decap ( \cM^\cU_{\ga + \eta} )$ into
$\cL_{\ga + \eta}$ that is cofinal:
$$\sup ( e_{\ga + \eta} \ '' \ \OR^{ \cM^\cU_{\ga + \eta} } )
= \OR^{ \cL_{\ga + \eta} }.$$\\

On the other hand, suppose that $\ga + \eta$ falls under Case B.
Then $e_{\ga + \eta}$ is a weak $\deg^\cU ( \ga + \eta )$-embedding
of $\cM^\cU_{\ga + \eta}$ into $\cL_{\ga + \eta}$.
Somewhat more elementarity will follow from the commutativity property
below.  For example, suppose that $\ga + \xi + 1$ is either
the last drop along $( \gb , \ga + \eta ]_U$, or
the least ordinal in $( \gb , \ga + \eta ]_U$
if no such drops exist.
Let $$Y(\cU , \ga + \eta) =
( i^\cU_{\ga + \xi +1 , \ga + \eta} \circ (i^*_{\ga + \xi + 1} )^\cU )
\, '' \, | (\cM^*_{\ga + \xi + 1} )^\cU | \ .$$
Then $e_{\ga + \eta}$ is $\gS_{deg(\ga + \eta ) + 1}$-elementary
on points in $Y ( \cU , \ga + \eta )$.\\}

\item[
\underline{Commutativity}]{\ \\

Suppose that 
$(\gb , \ga + \eta]_U \cap \cD^\cU  = \emptyset$.
Then 
both $i^\cU_{\gb , \ga + \eta}$
and $i^\cV_{\gb_\ell , \ga + \eta}$ are defined, and
$$e_{\ga + \eta} \circ i^\cU_{\gb , \ga + \eta}
= i^\cU_{\gb_\ell , \ga + \eta} \circ \gs_\ell \circ e_\gb
= i^\cV_{\gb_\ell , \ga + \eta} \circ \gs_\ell \res |\cR_\gb |$$\\

Now suppose that $\gb <_U (\ga + \xi ) <_U (\ga + \eta )$
and that $(\ga + \xi , \ga + \eta ]_U \cap \cD^\cU = \emptyset$.
Then both $i^U_{\ga + \xi , \ga + \eta}$
and $i^V_{\ga + \xi , \ga + \eta}$
are defined,
and 
$$e_{\ga + \eta} \circ i^U_{\ga + \xi , \ga + \eta}
= i^V_{\ga + \xi , \ga + \eta} \circ e_{\ga + \xi}.$$\\

Finally, suppose that $\ga + \eta \in \cD^\cU$ and that 
$\xi = \pred^\cU ( \ga + \eta )$.
Then 
$\xi = \pred^\cV ( \ga + \eta )$
and
$$e_{\ga + \eta} \circ ( i^*_{\ga + \eta} )^\cU
= ( i^*_{\ga + \eta} )^\cV \circ e_{\ga + \xi} \ .$$}

\end{list}

We now describe the construction of 
$\cV$ and the enlargement maps.  The last case is a bit 
trickier than the first three.

\vskip 3ex
\noindent
{\sc Case 1.}
Defining $\cV \res ( \ga + \eta)$ when $\ga + \eta$ is a limit ordinal.
\vskip 3ex

Here $\cV \res ( \ga + \eta )$
is just the limit of $\cV \res ( \ga + \xi)$ for $\ga + \xi < \ga + \eta$.

\vskip 3ex
\noindent
{\sc Case 2.}
Defining $\cV \res (\ga + \eta + 1)$ when
$\ga + \eta$ is a limit ordinal.
\vskip 3ex

Let $\gb = \rt^\cU ( \ga + \eta)$ and
$\gb_0 > \cdots \gb_\ell = \gb$, $\eta_0 < \cdots < \eta_\ell = \eta$
be the trace of $\eta$.
By a straightforward generalization of the Strong Uniqueness Theorem
6.2 of \cite{MiSt}, we see that 
$[\gb , \ga + \eta )_U$ is the sole cofinal wellfounded branch in
$\cU \res ( \ga + \eta )$.  The inductively maintained properties
of $\cV \res (\ga + \xi )$ and $e_{\ga + \xi}$ for $\ga + \xi
< \ga + \eta$ guarantee that 
$$\{ \  \ga + \gz \ | \mbox{ there is some } (\ga + \xi) <_U (\ga + \eta )
			\mbox{ s.t. }  (\ga + \gz) <_V ( \ga + \xi) \ \}$$
$$= [\gb_\ell , \ga + \eta_{\ell -1} +1 ]_V \ \cup \ 
(\ga + \eta_{\ell - 1} + 1 , \ga + \eta )_U$$
is the unique cofinal wellfounded branch in $\cV \res ( \ga + \eta )$.
This determines $\cV \res ( \ga + \eta )$,
$\cL_{\ga + \eta}$, and $e_{\ga + \eta}$ as limits in the usual way.

\vskip 3ex
\noindent
{\sc Case 3.}
Defining $\cV \res (\ga + \eta + 2)$
when either $\ga + \eta$ falls under Case B or
$E^\cU_{\ga + \eta} \not= \Fdot^{\cM^\cU_{\ga + \eta}}$.
\vskip 3ex

Let $\eta^* = \pred^\cU ( \ga + \eta + 1 )$.
Here we use the usual copying construction, using the maps
$e_{\eta^*}$ and $e_{\ga + \eta}$ to define $\cM^\cV_{\ga + \eta +1}$
as the model following $\cM^\cV_{\eta^*}$ in $\cV$.
If $E^\cU_{\ga + \eta} \not= \Fdot^{\cM^\cU_{\ga + \eta}}$,
then we put
$E^\cV_{\ga + \eta } = e_{\ga + \eta} ( E^\cU_{\ga + \eta} )$.
Suppose, on the other hand, that
$E^\cU_{\ga + \eta} = \Fdot^{\cM^\cU_{\ga + \eta}}$.
Then our case hypothesis implies that $\ga + \eta$ falls under Case B,
so $e_{\ga + \eta}$ is a weak $\deg^\cU ( \ga + \eta )$-embedding,
and more, of $\cM^\cU_{\ga + \eta}$ into $\cM^\cV_{\ga + \eta}$;
in particular, it has some elementarity with respect to the last
predicate of $\cM^\cU_{\ga + \eta}$.
Thus, it is reasonable that we set
$E^\cV_{\ga + \eta } = \Fdot^{\cM^\cV_{\ga + \eta}}$.

Having defined $E^\cV_{\ga + \eta }$, the 
rules for forming iteration trees determine $\cV \res ( \eta + 2 )$.
The map $e_{\ga + \eta + 1} : \cM^\cV_{\ga + \eta} \lra \cL_{\ga + \eta}$
is as given by the Shift Lemma 5.2 of \cite{MiSt}.  
It is relatively easy to check
that the inductive properties continue to hold at $\ga + \eta + 1$. 

\vskip 3ex
\noindent
{\sc Case 4.} Defining $\cV \res ( \ga + \eta + 2)$
when $\ga + \eta$ falls under Case A and 
$E^\cU_{\ga + \eta} = \Fdot^{\cM^\cU_{\ga + \eta}}$.
\vskip 3ex

We let $\gb_0 > \cdots > \gb_\ell$, $\eta_0 < \cdots < \eta_\ell$
be the trace of $\ga + \eta + 1$.
If $\gb = \pred^\cU (\ga + \eta + 1 )$,
then $\gb = \gb_0 < \ga$, $\ga + \eta = \ga + \eta_{\ell - 1}$,
and $\ga + \eta + 1 = \ga + \eta_\ell$.

Here we set $\cM^\cV_{ \ga + \eta + 1} = \cM^\cV_{\ga + \eta}$
and $\ga + \eta = \pred^\cV ( \ga + \eta + 1)$.
So it seems that we are merely padding from $\ga + \eta $
to $\ga + \eta + 1$ in the formation of $\cV$.
This is not exactly true,
for we also set $\gn^\cV_{ \ga + \eta} 
= e_{\ga + \eta} ( \gn^\cU_{\ga + \eta}  )$.
Technically, this is not padding in the sense of \cite{MiSt},
since we have given a condition by which it is possible that
$\ga + \eta = \pred^\cV (\ga + \gz )$ for some $\ga + \gz
> \ga + \eta + 1$; the condition
is that $\gn^\cV_{\ga + \xi} \leq \crit (E^\cU_{\ga + \gz})
<\gn_{\ga + \eta}^\cV$ 
for every $\ga + \xi < \ga + \eta $.
But clearly this is just a question of 
indexing; any phalanx iterable in the usual sense is iterable
in a sense involving this generalized form of padding.

Given this definition of $\cM^\cV_{\ga + \eta + 1}$,
we must next determine a map 
$$e_{\ga + \eta + 1} : | \cM^\cU_{\ga + \eta + 1} |
\lra | \cM^\cV_{\ga + \eta} | \ .$$
But first some motivation.
Consider the very simple situation in which $\eta_0 = 1$,
$\ell = 1$, and $m_{\gb_0}$, $m_{\gb_1}$, $n_{\gb_0}$,
and $n_{\gb_0}$ are all equal to $0$.  
Put $G = E^\cU_\ga$ and $F = \Fdot^{\cR_{\gb_0}}$.
Compare the sequence of extenders used along:
$$
| \cR_{\gb_0} |
\buildrel G
\over \lra
| \cM^\cU_{\ga + 1} |
\subset
| \cM^\cU_{\ga +2} | = | \ult ( \cR_\gb , i_G \, '' F ) |$$
with those used along:
$$ \cS_{\gb_1}
\buildrel F
\over \lra
\cS_{\gb_0}
\buildrel G
\over \lra
\cM^\cV_{\ga + 1}$$
The sequence of extenders used in the bottom iteration consists of images
of those used along the top system applied in a different order.
Note that 
$\crit (i_G \ '' F) = \crit (F) = \gp (\gk_{\gb_1})$.
Suppose that $x = i_{i_G \ '' F }(f)(a)$ for some 
$a \in [ \lh ( i_G \ '' F ) ]^{< \go}$ and $f \in | \cR_{\gb_1} |$
with $\dom (f) = [ \crit ( F ) ]^{|a|}$.
There are $b \in [ \lh (G) ]^{<\go}$ and $g \in |\cR_{\gb_1} |$
with $\dom ( g) = [\crit (G) ]^{|b|}$ such that 
$a = i_G (g)(b)$.
Then
$i_G ( i_F (f) \circ g )(b) \in |\cM^\cV_{\ga +1} |$.
This would be our value for 
$e_{\ga +2} (x)$. 

A bit more motivation for what is to come
can be gotten by considering Figure 1, which illustrates 
the situation when $\ell = 2$.
Compare the sequence of extenders that are used along

\begin{tabbing}
xxx\= xxxxxxxxxxx\= xxxxxxxxx\= \kill
\> $| \cR_{\gb_0} |
\lra
\buildrel i^\cU_{\gb_0 , \ga + \eta_0} \over \cdots
\lra
| \cM^\cU_{\ga + \eta_0} |$ \\
\ \\
\> \> $\subset
| \ \ult \, ( \, \cR_{\gb_1} \, , \, i^\cU_{\gb_0} {''} \Fdot^{\cR_{\gb_0}} \,
) \ | 
= | \cM^\cU_{\ga + \eta_0 + 1} |
\lra 
\buildrel i^\cU_{\ga + \eta_0 + 1 , \ga + \eta_1}
\over \cdots
\lra
| \cM^\cU_{\ga + \eta_1} |$ \\
\ \\
\> \> \> $\subset
| \, 
\ult \, (\, \cR_{\gb_2} \, , 
\, i^\cU_{\gb_1 , \ga + \eta_1} { '' } \Fdot^{\cR_{\gb_1}} \, ) \ |
= | \cM^\cU_{\ga + \eta_1 + 1} |$\\
\ \\
\end{tabbing}
with those used along
$$\cS_{\gb_2}
\buildrel \Fdot^{ \cR_{\gb_1} }
\over \ra
\cS_{\gb_1}
\buildrel \Fdot^{ \cR_{\gb_0} }
\over \ra
\cS_{\gb_0}
\ra
\buildrel i^\cV_{\gb_0 , \ga + \eta_0}
\over \cdots
\ra
\cM^\cV_{\ga + \eta_0}
=
\cM^\cV_{\ga + \eta_0 + 1}
\ra
\buildrel i^\cV_{\ga + \eta_0 + 1, \ga + \eta_1}
\over \cdots
\ra
\cM^\cV_{\ga + \eta_1}$$
Again, the second iteration is just the first
done in a different order, and then pushed over using the
maps $e_{\ga + \xi}$.

Now we give the general definition of $e_{\ga + \eta + 1}$ for Case 4.
Let $x \in | \cM^\cU_{\ga + \eta + 1} |$.
We may write $x$ as 
$$i^\cU_{\gb , \ga + \eta + 1} ( f) (a)$$
for some $a \in [\gn^\cU_{\ga + \eta} ]^{< \go}$
and function $f$ with $\dom (f ) = [ \gp ( \gk_\gb ) ]^{|a|}$.
If $m_\gb = 0$,
then we take $f \in | \cR_\gb |$.
If $m_\gb > 0$,
then $f : u \mapsto \gt_\psi^{\cR_\gb} [ u , q ]$
for some $\gS_{m_\gb}$-Skolem term $\gt_\psi$ and
parameter $q \in | \cR_\gb |$,
and by 
$$i^\cU_{\gb , \ga + \eta + 1} ( f) (a)$$
we really mean
$$\left[ a , f \right]^{\cR_\gb}_{E^\cU_{\ga + \eta}}
= \gt_\psi^{\cM^\cU_{\ga + \eta + 1}}
\left[ a , i^\cU_{\gb , \ga + \eta + 1} ( q) \right]$$
(we maintain this convention in what follows).

Next, consider any sequence with
$(f_\ell  , \dots , f_0 )$,
$( a_\ell , \dots , a_1 )$, 
and
$( b_{\ell - 1} , \dots , b_0 )$
such that:
\begin{list}{}{}
\item[1.]{$f_\ell = f$ and $a_\ell = a$.}
\item[2.]{If $\ell \geq 2$, then for $i = \ell -1 , \dots , 1$,
$$a_i \in [ \gn^\cU_{\ga + \eta_{i-1} }]^{<\go}$$
and 
$$b_i \in 
\left[ \sup \left( 
\{ \ \gn^\cU_\xi  \, | \, \xi \leq_U \ga + \eta_i \  \} \right)
- \gn^\cU_{\ga + \eta_{i-1}} 
\right]^{< \go} \ .$$
Also,
$$b_0 \in \left[ \sup \left( 
\{ \ \gn^\cU_\xi  \, | \, \xi \leq_U \ga + \eta_0 \  \} \right)
\right]^{< \go} \ .$$}
\item[3.]{For $i = \ell - 1 , \dots , 0$,
$f_i$ is a function, $f_i \in | \cR_{\gb_i} |$, and
$\dom ( f_i ) = [ \gp ( \gk_{ \gb_i } ) ]^{|a_i|}$.}
\item[4.]{If 
$\ell \geq 2$, then for 
$i = \ell , \dots , 2$,
$$a_i = i^\cU_{\gb_{i - 1}, \ga + \eta_{i - 1} } (f_{i -1} ) 
(a_{i-1} ) (b_{i-1} ) \ .$$}
\item[5.]{
$a_1 = i^\cU_{\gb_0 , \ga + \eta_0} ( f_0 )(b_0 )$}
\end{list}
When $i = \ell -1 , \dots , 0$ (as in 3 above)
we are able to choose $f_i \in | \cR_{\gb_i} |$,
since, as we already noted, $\deg^\cU (\ga + \eta_i ) = 0$.
Recall that we took
$\gs_i$ to be the embedding of $\cS_{\gb_i}$ into $\cS_{\gb_0}$.
Define
$$ g = ( \gs_\ell ( f_\ell ) ) \circ
	\dots
	\circ
	(\gs_1 ( f_1 )) \circ f_0$$
Note that 
if $\deg^\cV ( \ga + \eta ) = n_{\gb_0} = 0$,
then $g$ is an element of $\cS_{\gb_0}$;
otherwise,
$\deg^\cV ( \ga + \eta ) = n_{\gb_0} > 0$,
and $g$ is given by a $\gS_{n_{\gb_0}}$-Skolem term and parameter in
$\cS_{\gb_0}$.
For $i= \ell -1 , \dots 0$,
let
$$c_i
= 
e_{\ga + \eta_i } ( b_i )
=
e_{\ga + \eta_{\ell - 1} } (b_i) \ .$$
Finally, we set
$$e_{\ga + \eta + 1} ( x)
=
i^\cV_{
\gb_0, \ga + \eta} ( g ) (c_0 ) \cdots (c_{\ell -1}) \ .$$

The argument that $e_{\ga + \eta +1}$ is well-defined is straightforward.
Our rules determine that 
$$\cL_{\ga + \eta + 1} = i^\cV_{\gb_0 , \ga + \eta} (\gs_\ell (\cL_\gb ))$$
and a similar straightforward calculation shows that $e_{\ga + \eta + 1}$
has the required elementarity (almost by design).
The agreement property follows as usual, using the fact that if
$\xi < \lh (E^\cU_{\ga + \eta} )$,
then
$\xi = i^\cU_{\gb , \ga + \eta +1} (u \mapsto u_0 )(\{\xi \} )$.
The commutativity property also follows as usual, using the 
fact that $i^\cU_{\gb , \ga + \eta +1} ( x )
= i^\cU_{\gb , \ga + \eta +1} ( u \mapsto x )(a)$.
This completes the definition of our enlargement.

Now suppose, for contradiction, 
that $\cU$ has limit length, but that there is no
cofinal, wellfounded branch through $\cU$.
By (6)$_\ga$, there is a cofinal, wellfounded branch $b$ through $\cV$.
Let $\cG$ be the collection 
of $\ga + \eta + 1 \in b$
such that
$\ga + \eta $
falls under Case A and $E^\cU_{\ga + \eta} = \Fdot^{\cM^\cU_{\ga + \eta}}$
(i.e., the definition of $\cV \res (\ga + \eta + 2 )$ is by case 4).
Let us suppose that $\cG \not= \emptyset$.
If $\ga + \xi + 1 <_V \ga + \eta + 1$ are both in $\cG$,
then $\rt^\cU ( \ga + \xi + 1 ) > \rt^\cU ( \ga + \eta + 1 )$.
Therefore $\cG$ is finite;
let $\ga + \eta +1$ be the largest element of
$\cG$.  Suppose that the trace of $\ga + \eta + 1$ is
$\gb_0 > \cdots > \gb_\ell$,  $\eta_0 < \cdots < \eta_\ell$.
Then $\gb_0 = \min (b)$,
$$\cG = \{ \ \ga + \eta_{\ell - 1} + 1 , \dots , \ga + \eta_0 + 1 \ \} \ ,$$
and 
$$c = \{ \gb_0 \} \cup (b- (\ga + \eta +1))$$
is a cofinal branch through $\cU$.
Let $e$
be the direct limit of the maps $e_\xi$ for $\xi \in c$.
Then $e$ is an embedding of $\cM^\cU_c$ into $\cM^\cV_b$;
since the latter structure is wellfounded, both are.
This is a contradiction.  Similar arguments show that
$\cU$ is well-behaved in the remaining cases as well.
\qed {\tiny \ \ Lemma 3.18}
\vskip 3ex
 
\noindent
{\sc Lemma 3.19.}
Suppose that (1)$_\gb$ and (2)$_\gb$ holds for every $\gb < \ga$.
Then (6)$_\ga$ holds.
\vskip 3ex

\noindent
{\sc Proof.}
Consider any $\gb < \ga$. 
By Corollary 3.12, there is an iteration tree
$\cU_\gb$ on $W$ 
of length $\varphi_\gb + 1$
such that all extenders used on $\cU_\gb$ have length at least $\gL_\gb$,
and there is an elementary
embedding 
$j_\gb : \cS_\gb \lra \cM^{\cU_\gb}_{\varphi_\gb}$
with critical point at least $\gp ( \gk_\gb )$.
It follows that $(( \langle \ \cM^{\cU_\gb}_{\varphi_\gb} \ |
\ \gb < \ga \ \rangle , W ), \gLvec \res \ga )$
is a phalanx.
The proof of the following fact is 
is similar to the proof of the main result in Section 9 of \cite{St1};
we omit the details.

\vskip 3ex
\noindent
{\sc Fact 3.19.1.}
$(( \langle \ \cM^{\cU_\gb}_{\varphi_\gb} \ |
\ \gb < \ga \ \rangle , W ) , \gLvec \res \ga )$
is an iterable phalanx.
\vskip 3ex

The system of maps
$(\langle \ j_\gb \ | \  \gb < \ga \ \rangle , \id \res W )$
now guarantees that 
$( (\cSvec \res \ga , W ) , \gLvec \res \ga )$ is 
also iterable in the usual way.
\qed {\tiny \ \ Lemma 3.19}
\vskip 3ex

This completes the proof that the implications listed before
Lemma 3.9 hold, and consequently that (1)$_\ga$ through (6)$_\ga$
hold for every $\ga < \gg$.  As already explained (just after the 
statement of Corollary 3.12), Theorem 1.1 follows.
\qed {\tiny \ \ Theorem 1.1}

\vskip 4ex

\noindent
{\small {\sc Department of Mathematics, University of Florida, 
Gainesville, FL \ 32611\\
Department of Mathematics, MIT, 
Cambridge, MA \ 02139\\
Department of Mathematics, UCLA,
Los Angeles, CA \ 90024\\}}

\noindent
{\small {\it E-mail addresses:} \ mitchell@@math.ufl.edu, 
ernest@@math.mit.edu, steel@@math.ucla.edu}
\end{document}